
\magnification=1200
\hsize=13.50cm
\vsize=18cm
\hoffset=-2mm
\voffset=.8cm
\parindent=12pt   \parskip=0pt
\hfuzz=1pt

\pretolerance=500 \tolerance=1000  \brokenpenalty=5000

\catcode`\@=11

\font\eightrm=cmr8
\font\eighti=cmmi8
\font\eightsy=cmsy8
\font\eightbf=cmbx8
\font\eighttt=cmtt8
\font\eightit=cmti8
\font\eightsl=cmsl8
\font\sevenrm=cmr7
\font\seveni=cmmi7
\font\sevensy=cmsy7
\font\sevenbf=cmbx7

\font\sixrm=cmr6
\font\sixi=cmmi6
\font\sixsy=cmsy6
\font\sixbf=cmbx6

\font\douzebf=cmbx10 at 12pt

\font\twelvebf=cmbx10 at 12pt

\font\tencal=eusm10

\font\sevencal=eusm7

\font\fivecal=eusm5
\newfam\calfam
\textfont\calfam=\tencal
\scriptfont\calfam=\sevencal
\scriptscriptfont\calfam=\fivecal
\def\cal#1{{\fam\calfam\relax#1}}

\skewchar\eighti='177 \skewchar\sixi='177
\skewchar\eightsy='60 \skewchar\sixsy='60

\def\tenpoint{%
   \textfont0=\tenrm \scriptfont0=\sevenrm
   \scriptscriptfont0=\fiverm
   \def\rm{\fam\z@\tenrm}%
   \textfont1=\teni  \scriptfont1=\seveni
   \scriptscriptfont1=\fivei
   \def\oldstyle{\fam\@ne\teni}\let\old=\oldstyle
   \textfont2=\tensy \scriptfont2=\sevensy
   \scriptscriptfont2=\fivesy
   \textfont\itfam=\tenit
   \def\it{\fam\itfam\tenit}%
   \textfont\slfam=\tensl
   \def\sl{\fam\slfam\tensl}%
   \textfont\bffam=\tenbf
   \scriptfont\bffam=\sevenbf
   \scriptscriptfont\bffam=\fivebf
   \def\bf{\fam\bffam\tenbf}%
   \textfont\ttfam=\tentt
   \def\tt{\fam\ttfam\tentt}%
   \abovedisplayskip=12pt plus 3pt minus 9pt
   \belowdisplayskip=\abovedisplayskip
   \abovedisplayshortskip=0pt plus 3pt
   \belowdisplayshortskip=4pt plus 3pt
   \smallskipamount=3pt plus 1pt minus 1pt
   \medskipamount=6pt plus 2pt minus 2pt
   \bigskipamount=12pt plus 4pt minus 4pt
   \normalbaselineskip=12pt
   \setbox\strutbox=\hbox{\vrule height8.5pt depth3.5pt width0pt}%
   \let\bigf@nt=\tenrm
   \let\smallf@nt=\sevenrm
   \normalbaselines\rm}

\def\eightpoint{%
   \textfont0=\eightrm \scriptfont0=\sixrm
   \scriptscriptfont0=\fiverm
   \def\rm{\fam\z@\eightrm}%
   \textfont1=\eighti  \scriptfont1=\sixi
   \scriptscriptfont1=\fivei
   \def\oldstyle{\fam\@ne\eighti}\let\old=\oldstyle
   \textfont2=\eightsy \scriptfont2=\sixsy
   \scriptscriptfont2=\fivesy
   \textfont\itfam=\eightit
   \def\it{\fam\itfam\eightit}%
   \textfont\slfam=\eightsl
   \def\sl{\fam\slfam\eightsl}%
   \textfont\bffam=\eightbf
   \scriptfont\bffam=\sixbf
   \scriptscriptfont\bffam=\fivebf
   \def\bf{\fam\bffam\eightbf}%
   \textfont\ttfam=\eighttt
   \def\tt{\fam\ttfam\eighttt}%
   \abovedisplayskip=9pt plus 3pt minus 9pt
   \belowdisplayskip=\abovedisplayskip
   \abovedisplayshortskip=0pt plus 3pt
   \belowdisplayshortskip=3pt plus 3pt
   \smallskipamount=2pt plus 1pt minus 1pt
   \medskipamount=4pt plus 2pt minus 1pt
   \bigskipamount=9pt plus 3pt minus 3pt
   \normalbaselineskip=9pt
   \setbox\strutbox=\hbox{\vrule height7pt depth2pt width0pt}%
   \let\bigf@nt=\eightrm
   \let\smallf@nt=\sixrm
   \normalbaselines\rm}
\tenpoint


\def\pc#1{\bigf@nt#1\smallf@nt}

\catcode`\;=\active
\def;{\relax\ifhmode\ifdim\lastskip>\z@\unskip\fi \kern\fontdimen2
 -1.2 \fontdimen3 \string;}

\catcode`\:=\active
\def:{\relax\ifhmode\ifdim\lastskip>\z@\unskip\fi\penalty\@M\
\fi\string:}

\catcode`\!=\active
\def!{\relax\ifhmode\ifdim\lastskip>\z@ \unskip\fi\kern\fontdimen2
 -1.1 \fontdimen3 \string!}

\catcode`\?=\active
\def?{\relax\ifhmode\ifdim\lastskip>\z@ \unskip\fi\kern\fontdimen2
 -1.1 \fontdimen3 \string?}

\frenchspacing

\newtoks\auteurcourant
\auteurcourant={\hfil}

\newtoks\titrecourant
\titrecourant={\hfil}

\newtoks\hautpagetitre
\hautpagetitre={\hfil}

\newtoks\baspagetitre
\baspagetitre={\hfil}

\newtoks\hautpagegauche
\hautpagegauche={\eightpoint\rlap{\folio}\hfil\the
\auteurcourant\hfil}

\newtoks\hautpagedroite
\hautpagedroite={\eightpoint\hfil\the
\titrecourant\hfil\llap{\folio}}

\newtoks\baspagegauche
\baspagegauche={\hfil}

\newtoks\baspagedroite
\baspagedroite={\hfil}

\newif\ifpagetitre
\pagetitretrue


\headline={\ifpagetitre\the\hautpagetitre
\else\ifodd\pageno\the\hautpagedroite\else\the
\hautpagegauche\fi\fi}

\footline={\ifpagetitre\the\baspagetitre\else
\ifodd\pageno\the\baspagedroite\else\the
\baspagegauche\fi\fi
\global\pagetitrefalse}

\def\raggedbottom{\topskip 10pt plus 36pt\r@ggedbottomtrue}

\def\point{\raise.2ex\hbox{\douzebf .}}
\def\pointir{\unskip . --- \ignorespaces}


\def\Medbreak{\vskip-\lastskip\medbreak}

\def\rem#1\endrem{%
\Medbreak {\it#1\unskip} : }


\long\def\th#1 #2\enonce#3\endth{%
\Medbreak {\pc#1} {#2\unskip}\pointir{\it #3}\medskip}


\long\def\thm#1 #2\enonce#3\endthm{%
\Medbreak {\pc#1} {#2\unskip}\pointir{\it #3}\medskip}

\def\decale#1{\smallbreak\hskip 28pt\llap{#1}\kern 5pt}

\def\decaledecale#1{\smallbreak\hskip 34pt\llap{#1}\kern 5pt}

\let\@ldmessage=\message

\def\message#1{{\def\pc{\string\pc\space}%
\def\'{\string'}\def\`{\string`}%
\def\^{\string^}\def\"{\string"}%
\@ldmessage{#1}}}

\def\({{\rm (}}

\def\){{\rm )}}


\def\up#1{\raise 1ex\hbox{\smallf@nt#1}}

\def\diagram#1{\def\normalbaselines{\baselineskip=0pt
\lineskip=5pt}\matrix{#1}}


\def\longmaprightover#1#2{\smash{\mathop{\hbox to#2
{\rightarrowfill}}\limits^{\scriptstyle#1}}}

\def\longmapleftover#1#2{\smash{\mathop{\hbox to#2
{\leftarrowfill}}\limits^{\scriptstyle#1}}}

\def\longmaprightunder#1#2{\smash{\mathop{\hbox to#2
{\rightarrowfill}}\limits_{\scriptstyle#1}}}

\def\longmapleftunder#1#2{\smash{\mathop{\hbox to#2
{\leftarrowfill}}\limits_{\scriptstyle#1}}}

\def\longhookrightarrowover#1#2{\smash{\mathop{\lhook\joinrel
\mathrel{\hbox to #2{\rightarrowfill}}}\limits^{\scriptstyle#1}}}

\def\longhookrightarrowunder#1#2{\smash{\mathop{\lhook\joinrel
\mathrel{\hbox to #2{\rightarrowfill}}}\limits_{\scriptstyle#1}}}

\def\longhookleftarrowover#1#2{\smash{\mathop{{\hbox to #2
{\leftarrowfill}}\joinrel\kern -0.9mm\mathrel\rhook}
\limits^{\scriptstyle#1}}}

\def\longhookleftarrowunder#1#2{\smash{\mathop{{\hbox to #2
{\leftarrowfill}}\joinrel\kern -0.9mm\mathrel\rhook}
\limits_{\scriptstyle#1}}}

\def\longtwoheadrightarrowover#1#2{\smash{\mathop{{\hbox to #2
{\rightarrowfill}}\kern -3.25mm\joinrel\mathrel\rightarrow}
\limits^{\scriptstyle#1}}}

\def\longtwoheadrightarrowunder#1#2{\smash{\mathop{{\hbox to #2
{\rightarrowfill}}\kern -3.25mm\joinrel\mathrel\rightarrow}
\limits_{\scriptstyle#1}}}

\def\longtwoheadleftarrowover#1#2{\smash{\mathop{\joinrel\mathrel
\leftarrow\kern -3.8mm{\hbox to #2{\leftarrowfill}}}
\limits^{\scriptstyle#1}}}

\def\longtwoheadleftarrowunder#1#2{\smash{\mathop{\joinrel\mathrel
\leftarrow\kern -3.8mm{\hbox to #2{\leftarrowfill}}}
\limits_{\scriptstyle#1}}}


\def\longmapsto#1{\mapstochar\mathrel{\joinrel \kern-0.2mm\hbox to
#1mm{\rightarrowfill}}}

\def\og{\leavevmode\raise.3ex\hbox{$\scriptscriptstyle
\langle\!\langle\,$}}
\def\fg{\leavevmode\raise.3ex\hbox{$\scriptscriptstyle
\,\rangle\!\rangle$}}

\def\section#1#2{\vskip 5mm {\bf {#1}. {#2}}\vskip 5mm}
\def\subsection#1#2{\vskip 3mm {\it #2}\vskip 3mm}

\catcode`\@=12

\showboxbreadth=-1  \showboxdepth=-1


\message{`lline' & `vector' macros from LaTeX}

\def\Grille{\setbox13=\vbox to 5\unitlength{\hrule width 109mm \vfill}
\setbox13=\vbox to 65mm
{\offinterlineskip\leaders\copy13\vfill\kern-1pt\hrule}
\setbox14=\hbox to 5\unitlength{\vrule height 65mm\hfill}
\setbox14=\hbox to 109mm{\leaders\copy14\hfill\kern-2mm \vrule height
65mm}
\ht14=0pt\dp14=0pt\wd14=0pt \setbox13=\vbox to 0pt
{\vss\box13\offinterlineskip\box14} \wd13=0pt\box13}

\def\rule(#1,#2)\dir(#3,#4)\long#5{%
\noalign{\leftput(#1,#2){\lline(#3,#4){#5}}}}

\def\arrow(#1,#2)\dir(#3,#4)\length#5{%
\noalign{\leftput(#1,#2){\vector(#3,#4){#5}}}}

\def\put(#1,#2)#3{\noalign{\setbox1=\hbox{%
\kern #1\unitlength \raise #2\unitlength\hbox{$#3$}}%
\ht1=0pt \wd1=0pt \dp1=0pt\box1}}

\catcode`@=11

\def\{{\relax\ifmmode\lbrace\else$\lbrace$\fi}

\def\}{\relax\ifmmode\rbrace\else$\rbrace$\fi}

\def\newcount{\alloc@0\count\countdef\insc@unt}

\def\newdimen{\alloc@1\dimen\dimendef\insc@unt}

\def\newwrite{\alloc@7\write\chardef\sixt@@n}

\newwrite\@unused
\newcount\@tempcnta
\newcount\@tempcntb
\newdimen\@tempdima
\newdimen\@tempdimb
\newbox\@tempboxa

\def\@spaces{\space\space\space\space}

\def\@whilenoop#1{}

\def\@whiledim#1\do #2{\ifdim #1\relax#2\@iwhiledim{#1\relax#2}\fi}

\def\@iwhiledim#1{\ifdim #1\let\@nextwhile=\@iwhiledim
\else\let\@nextwhile=\@whilenoop\fi\@nextwhile{#1}}

\def\@badlinearg{\@latexerr{Bad \string\line\space or \string\vector
\space argument}}

\def\@latexerr#1#2{\begingroup
\edef\@tempc{#2}\expandafter\errhelp\expandafter{\@tempc}%

\def\@eha{Your command was ignored.^^JType \space I <command> <return>
\space to replace it with another command,^^Jor \space <return> \space
to continue without it.}

\def\@ehb{You've lost some text. \space \@ehc}

\def\@ehc{Try typing \space <return> \space to proceed.^^JIf that
doesn't work, type \space X <return> \space to quit.}

\def\@ehd{You're in trouble here.  \space\@ehc}

\typeout{LaTeX error.  \space See LaTeX manual for explanation.^^J
\space\@spaces\@spaces\@spaces Type \space H <return> \space for
immediate help.}\errmessage{#1}\endgroup}

\def\typeout#1{{\let\protect\string\immediate\write\@unused{#1}}}

\font\tenln = line10
\font\tenlnw = linew10

\newdimen\@wholewidth
\newdimen\@halfwidth
\newdimen\unitlength

\unitlength =1pt

\def\thinlines{\let\@linefnt\tenln \let\@circlefnt\tencirc
\@wholewidth\fontdimen8\tenln \@halfwidth .5\@wholewidth}

\def\thicklines{\let\@linefnt\tenlnw \let\@circlefnt\tencircw
\@wholewidth\fontdimen8\tenlnw \@halfwidth .5\@wholewidth}

\def\linethickness#1{\@wholewidth #1\relax \@halfwidth .5
\@wholewidth}

\newif\if@negarg

\def\lline(#1,#2)#3{\@xarg #1\relax \@yarg #2\relax
\@linelen=#3\unitlength \ifnum\@xarg =0 \@vline \else \ifnum\@yarg =0
\@hline \else \@sline\fi \fi}

\def\@sline{\ifnum\@xarg< 0 \@negargtrue \@xarg -\@xarg \@yyarg
-\@yarg \else \@negargfalse \@yyarg \@yarg \fi
\ifnum \@yyarg >0 \@tempcnta\@yyarg \else \@tempcnta - \@yyarg \fi
\ifnum\@tempcnta>6 \@badlinearg\@tempcnta0 \fi
\setbox\@linechar\hbox{\@linefnt\@getlinechar(\@xarg,\@yyarg)}%
\ifnum \@yarg >0 \let\@upordown\raise \@clnht\z@
\else\let\@upordown\lower \@clnht \ht\@linechar\fi
\@clnwd=\wd\@linechar
\if@negarg \hskip -\wd\@linechar \def\@tempa{\hskip -2\wd \@linechar}
\else \let\@tempa\relax \fi
\@whiledim \@clnwd <\@linelen \do {\@upordown\@clnht\copy\@linechar
\@tempa \advance\@clnht \ht\@linechar \advance\@clnwd \wd\@linechar}%
\advance\@clnht -\ht\@linechar \advance\@clnwd -\wd\@linechar
\@tempdima\@linelen\advance\@tempdima -\@clnwd
\@tempdimb\@tempdima\advance\@tempdimb -\wd\@linechar
\if@negarg \hskip -\@tempdimb \else \hskip \@tempdimb \fi
\multiply\@tempdima \@m\@tempcnta \@tempdima \@tempdima \wd\@linechar
\divide\@tempcnta \@tempdima \@tempdima \ht\@linechar
\multiply\@tempdima \@tempcnta \divide\@tempdima \@m \advance\@clnht
\@tempdima
\ifdim \@linelen <\wd\@linechar \hskip \wd\@linechar
\else\@upordown\@clnht\copy\@linechar\fi}

\def\@hline{\ifnum \@xarg <0 \hskip -\@linelen \fi
\vrule height \@halfwidth depth \@halfwidth width \@linelen
\ifnum \@xarg <0 \hskip -\@linelen \fi}

\def\@getlinechar(#1,#2){\@tempcnta#1\relax
\multiply\@tempcnta 8\advance\@tempcnta -9
\ifnum #2>0 \advance\@tempcnta #2\relax
  \else\advance\@tempcnta -#2\relax\advance\@tempcnta 64 \fi
\char\@tempcnta}

\def\vector(#1,#2)#3{\@xarg #1\relax \@yarg #2\relax
\@linelen=#3\unitlength
\ifnum\@xarg =0 \@vvector \else \ifnum\@yarg =0 \@hvector \else
\@svector\fi \fi}

\def\@hvector{\@hline\hbox to 0pt{\@linefnt \ifnum \@xarg <0
\@getlarrow(1,0)\hss\else \hss\@getrarrow(1,0)\fi}}

\def\@vvector{\ifnum \@yarg <0 \@downvector \else \@upvector \fi}

\def\@svector{\@sline\@tempcnta\@yarg \ifnum\@tempcnta <0
\@tempcnta=-\@tempcnta\fi \ifnum\@tempcnta <5 \hskip -\wd\@linechar
\@upordown\@clnht \hbox{\@linefnt \if@negarg
\@getlarrow(\@xarg,\@yyarg) \else \@getrarrow(\@xarg,\@yyarg)
\fi}\else\@badlinearg\fi}

\def\@getlarrow(#1,#2){\ifnum #2 =\z@ \@tempcnta='33\else
\@tempcnta=#1\relax\multiply\@tempcnta \sixt@@n \advance\@tempcnta -9
\@tempcntb=#2\relax \multiply\@tempcntb \tw@ \ifnum \@tempcntb >0
\advance\@tempcnta \@tempcntb\relax \else\advance\@tempcnta
-\@tempcntb\advance\@tempcnta 64 \fi\fi \char\@tempcnta}

\def\@getrarrow(#1,#2){\@tempcntb=#2\relax \ifnum\@tempcntb < 0
\@tempcntb=-\@tempcntb\relax\fi \ifcase \@tempcntb\relax
\@tempcnta='55 \or \ifnum #1<3 \@tempcnta=#1\relax\multiply\@tempcnta
24 \advance\@tempcnta -6 \else \ifnum #1=3 \@tempcnta=49
\else\@tempcnta=58 \fi\fi\or \ifnum #1<3
\@tempcnta=#1\relax\multiply\@tempcnta 24 \advance\@tempcnta -3 \else
\@tempcnta=51\fi\or \@tempcnta=#1\relax\multiply\@tempcnta \sixt@@n
\advance\@tempcnta -\tw@ \else \@tempcnta=#1\relax\multiply\@tempcnta
\sixt@@n \advance\@tempcnta 7 \fi \ifnum #2<0 \advance\@tempcnta 64
\fi \char\@tempcnta}

\def\@vline{\ifnum \@yarg <0 \@downline \else \@upline\fi}

\def\@upline{\hbox to \z@{\hskip -\@halfwidth \vrule width
\@wholewidth height \@linelen depth \z@\hss}}

\def\@downline{\hbox to \z@{\hskip -\@halfwidth \vrule width
\@wholewidth height \z@ depth \@linelen \hss}}

\def\@upvector{\@upline\setbox\@tempboxa
\hbox{\@linefnt\char'66}\raise \@linelen \hbox to\z@{\lower
\ht\@tempboxa \box\@tempboxa\hss}}

\def\@downvector{\@downline\lower \@linelen \hbox to
\z@{\@linefnt\char'77\hss}}

\thinlines

\newcount\@xarg
\newcount\@yarg
\newcount\@yyarg
\newcount\@multicnt
\newdimen\@xdim
\newdimen\@ydim
\newbox\@linechar
\newdimen\@linelen
\newdimen\@clnwd
\newdimen\@clnht
\newdimen\@dashdim
\newbox\@dashbox
\newcount\@dashcnt
\catcode`@=12

\newbox\tbox
\newbox\tboxa

\def\leftzer#1{\setbox\tbox=\hbox to 0pt{#1\hss}%
\ht\tbox=0pt \dp\tbox=0pt \box\tbox}

\def\rightzer#1{\setbox\tbox=\hbox to 0pt{\hss #1}%
\ht\tbox=0pt \dp\tbox=0pt \box\tbox}

\def\centerzer#1{\setbox\tbox=\hbox to 0pt{\hss #1\hss}%
\ht\tbox=0pt \dp\tbox=0pt \box\tbox}

\def\leftput(#1,#2)#3{\setbox\tboxa=\hbox{%
\kern #1\unitlength \raise #2\unitlength\hbox{\leftzer{#3}}}%
\ht\tboxa=0pt \wd\tboxa=0pt \dp\tboxa=0pt\box\tboxa}

\def\rightput(#1,#2)#3{\setbox\tboxa=\hbox{%
\kern #1\unitlength \raise #2\unitlength\hbox{\rightzer{#3}}}%
\ht\tboxa=0pt \wd\tboxa=0pt \dp\tboxa=0pt\box\tboxa}

\def\centerput(#1,#2)#3{\setbox\tboxa=\hbox{%
\kern #1\unitlength \raise #2\unitlength\hbox{\centerzer{#3}}}%
\ht\tboxa=0pt \wd\tboxa=0pt \dp\tboxa=0pt\box\tboxa}

\unitlength=1mm

\expandafter\ifx\csname amssym.def\endcsname\relax \else
\endinput\fi
%
\expandafter\edef\csname amssym.def\endcsname{%
        \catcode`\noexpand\@=\the\catcode`\@\space}
\catcode`\@=11
%

\def\undefine#1{\let#1\undefined}
\def\newsymbol#1#2#3#4#5{\let\next@\relax
  \ifnum#2=\@ne\let\next@\msafam@\else
  \ifnum#2=\tw@\let\next@\msbfam@\fi\fi
  \mathchardef#1="#3\next@#4#5}
\def\mathhexbox@#1#2#3{\relax
  \ifmmode\mathpalette{}{\m@th\mathchar"#1#2#3}%
  \else\leavevmode\hbox{$\m@th\mathchar"#1#2#3$}\fi}
\def\hexnumber@#1{\ifcase#1 0\or 1\or 2\or 3
\or 4\or 5\or 6\or 7\or 8\or
  9\or A\or B\or C\or D\or E\or F\fi}

\font\tenmsa=msam10
\font\sevenmsa=msam7
\font\fivemsa=msam5
\newfam\msafam
\textfont\msafam=\tenmsa
\scriptfont\msafam=\sevenmsa
\scriptscriptfont\msafam=\fivemsa
\edef\msafam@{\hexnumber@\msafam}
\mathchardef\dabar@"0\msafam@39
\def\dashrightarrow{\mathrel{\dabar@\dabar@\mathchar"0
\msafam@4B}}
\def\dashleftarrow{\mathrel{\mathchar"0\msafam@4C
\dabar@\dabar@}}

\def\ulcorner{\delimiter"4\msafam@70\msafam@70 }
\def\urcorner{\delimiter"5\msafam@71\msafam@71 }
\def\llcorner{\delimiter"4\msafam@78\msafam@78 }
\def\lrcorner{\delimiter"5\msafam@79\msafam@79 }
\def\yen{{\mathhexbox@\msafam@55}}
\def\checkmark{{\mathhexbox@\msafam@58}}
\def\circledR{{\mathhexbox@\msafam@72}}
\def\maltese{{\mathhexbox@\msafam@7A}}

\font\tenmsb=msbm10
\font\sevenmsb=msbm7
\font\fivemsb=msbm5
\newfam\msbfam
\textfont\msbfam=\tenmsb
\scriptfont\msbfam=\sevenmsb
\scriptscriptfont\msbfam=\fivemsb
\edef\msbfam@{\hexnumber@\msbfam}
\def\Bbb#1{{\fam\msbfam\relax#1}}
\def\widehat#1{\setbox\z@\hbox{$\m@th#1$}%
  \ifdim\wd\z@>\tw@ em\mathaccent"0\msbfam@5B{#1}%
  \else\mathaccent"0362{#1}\fi}
\def\widetilde#1{\setbox\z@\hbox{$\m@th#1$}%
  \ifdim\wd\z@>\tw@ em\mathaccent"0\msbfam@5D{#1}%
  \else\mathaccent"0365{#1}\fi}
\font\teneufm=eufm10
\font\seveneufm=eufm7
\font\fiveeufm=eufm5
\newfam\eufmfam
\textfont\eufmfam=\teneufm
\scriptfont\eufmfam=\seveneufm
\scriptscriptfont\eufmfam=\fiveeufm
\def\frak#1{{\fam\eufmfam\relax#1}}

\csname amssym.def\endcsname

\expandafter\ifx\csname pre amssym.tex at\endcsname\relax \else
\endinput\fi
\expandafter\chardef\csname pre amssym.tex at\endcsname=\the
\catcode`\@
\catcode`\@=11
\begingroup\ifx\undefined\newsymbol \else\def\input#1
{\endgroup}\fi
\input amssym.def \relax
\newsymbol\boxdot 1200
\newsymbol\boxplus 1201
\newsymbol\boxtimes 1202
\newsymbol\square 1003
\newsymbol\blacksquare 1004
\newsymbol\centerdot 1205
\newsymbol\lozenge 1006
\newsymbol\blacklozenge 1007
\newsymbol\circlearrowright 1308
\newsymbol\circlearrowleft 1309
\undefine\rightleftharpoons
\newsymbol\rightleftharpoons 130A
\newsymbol\leftrightharpoons 130B
\newsymbol\boxminus 120C
\newsymbol\Vdash 130D
\newsymbol\Vvdash 130E
\newsymbol\vDash 130F
\newsymbol\twoheadrightarrow 1310
\newsymbol\twoheadleftarrow 1311
\newsymbol\leftleftarrows 1312
\newsymbol\rightrightarrows 1313
\newsymbol\upuparrows 1314
\newsymbol\downdownarrows 1315
\newsymbol\upharpoonright 1316
  
\newsymbol\downharpoonright 1317
\newsymbol\upharpoonleft 1318
\newsymbol\downharpoonleft 1319
\newsymbol\rightarrowtail 131A
\newsymbol\leftarrowtail 131B
\newsymbol\leftrightarrows 131C
\newsymbol\rightleftarrows 131D
\newsymbol\Lsh 131E
\newsymbol\Rsh 131F
\newsymbol\rightsquigarrow 1320
\newsymbol\leftrightsquigarrow 1321
\newsymbol\looparrowleft 1322
\newsymbol\looparrowright 1323
\newsymbol\circeq 1324
\newsymbol\succsim 1325
\newsymbol\gtrsim 1326
\newsymbol\gtrapprox 1327
\newsymbol\multimap 1328
\newsymbol\therefore 1329
\newsymbol\because 132A
\newsymbol\doteqdot 132B
  
\newsymbol\triangleq 132C
\newsymbol\precsim 132D
\newsymbol\lesssim 132E
\newsymbol\lessapprox 132F
\newsymbol\eqslantless 1330
\newsymbol\eqslantgtr 1331
\newsymbol\curlyeqprec 1332
\newsymbol\curlyeqsucc 1333
\newsymbol\preccurlyeq 1334
\newsymbol\leqq 1335
\newsymbol\leqslant 1336
\newsymbol\lessgtr 1337
\newsymbol\backprime 1038
\newsymbol\risingdotseq 133A
\newsymbol\fallingdotseq 133B
\newsymbol\succcurlyeq 133C
\newsymbol\geqq 133D
\newsymbol\geqslant 133E
\newsymbol\gtrless 133F
\newsymbol\sqsubset 1340
\newsymbol\sqsupset 1341
\newsymbol\vartriangleright 1342
\newsymbol\vartriangleleft 1343
\newsymbol\trianglerighteq 1344
\newsymbol\trianglelefteq 1345
\newsymbol\bigstar 1046
\newsymbol\between 1347
\newsymbol\blacktriangledown 1048
\newsymbol\blacktriangleright 1349
\newsymbol\blacktriangleleft 134A
\newsymbol\vartriangle 134D
\newsymbol\blacktriangle 104E
\newsymbol\triangledown 104F
\newsymbol\eqcirc 1350
\newsymbol\lesseqgtr 1351
\newsymbol\gtreqless 1352
\newsymbol\lesseqqgtr 1353
\newsymbol\gtreqqless 1354
\newsymbol\Rrightarrow 1356
\newsymbol\Lleftarrow 1357
\newsymbol\veebar 1259
\newsymbol\barwedge 125A
\newsymbol\doublebarwedge 125B
\undefine\angle
\newsymbol\angle 105C
\newsymbol\measuredangle 105D
\newsymbol\sphericalangle 105E
\newsymbol\varpropto 135F
\newsymbol\smallsmile 1360
\newsymbol\smallfrown 1361
\newsymbol\Subset 1362
\newsymbol\Supset 1363
\newsymbol\Cup 1264
  
\newsymbol\Cap 1265
  
\newsymbol\curlywedge 1266
\newsymbol\curlyvee 1267
\newsymbol\leftthreetimes 1268
\newsymbol\rightthreetimes 1269
\newsymbol\subseteqq 136A
\newsymbol\supseteqq 136B
\newsymbol\bumpeq 136C
\newsymbol\Bumpeq 136D
\newsymbol\lll 136E
  
\newsymbol\ggg 136F
  
\newsymbol\circledS 1073
\newsymbol\pitchfork 1374
\newsymbol\dotplus 1275
\newsymbol\backsim 1376
\newsymbol\backsimeq 1377
\newsymbol\complement 107B
\newsymbol\intercal 127C
\newsymbol\circledcirc 127D
\newsymbol\circledast 127E
\newsymbol\circleddash 127F
\newsymbol\lvertneqq 2300
\newsymbol\gvertneqq 2301
\newsymbol\nleq 2302
\newsymbol\ngeq 2303
\newsymbol\nless 2304
\newsymbol\ngtr 2305
\newsymbol\nprec 2306
\newsymbol\nsucc 2307
\newsymbol\lneqq 2308
\newsymbol\gneqq 2309
\newsymbol\nleqslant 230A
\newsymbol\ngeqslant 230B
\newsymbol\lneq 230C
\newsymbol\gneq 230D
\newsymbol\npreceq 230E
\newsymbol\nsucceq 230F
\newsymbol\precnsim 2310
\newsymbol\succnsim 2311
\newsymbol\lnsim 2312
\newsymbol\gnsim 2313
\newsymbol\nleqq 2314
\newsymbol\ngeqq 2315
\newsymbol\precneqq 2316
\newsymbol\succneqq 2317
\newsymbol\precnapprox 2318
\newsymbol\succnapprox 2319
\newsymbol\lnapprox 231A
\newsymbol\gnapprox 231B
\newsymbol\nsim 231C
\newsymbol\ncong 231D
\newsymbol\diagup 201E
\newsymbol\diagdown 201F
\newsymbol\varsubsetneq 2320
\newsymbol\varsupsetneq 2321
\newsymbol\nsubseteqq 2322
\newsymbol\nsupseteqq 2323
\newsymbol\subsetneqq 2324
\newsymbol\supsetneqq 2325
\newsymbol\varsubsetneqq 2326
\newsymbol\varsupsetneqq 2327
\newsymbol\subsetneq 2328
\newsymbol\supsetneq 2329
\newsymbol\nsubseteq 232A
\newsymbol\nsupseteq 232B
\newsymbol\nparallel 232C
\newsymbol\nmid 232D
\newsymbol\nshortmid 232E
\newsymbol\nshortparallel 232F
\newsymbol\nvdash 2330
\newsymbol\nVdash 2331
\newsymbol\nvDash 2332
\newsymbol\nVDash 2333
\newsymbol\ntrianglerighteq 2334
\newsymbol\ntrianglelefteq 2335
\newsymbol\ntriangleleft 2336
\newsymbol\ntriangleright 2337
\newsymbol\nleftarrow 2338
\newsymbol\nrightarrow 2339
\newsymbol\nLeftarrow 233A
\newsymbol\nRightarrow 233B
\newsymbol\nLeftrightarrow 233C
\newsymbol\nleftrightarrow 233D
\newsymbol\divideontimes 223E
\newsymbol\varnothing 203F
\newsymbol\nexists 2040
\newsymbol\Finv 2060
\newsymbol\Game 2061
\newsymbol\mho 2066
\newsymbol\eth 2067
\newsymbol\eqsim 2368
\newsymbol\beth 2069
\newsymbol\gimel 206A
\newsymbol\daleth 206B
\newsymbol\lessdot 236C
\newsymbol\gtrdot 236D
\newsymbol\ltimes 226E
\newsymbol\rtimes 226F
\newsymbol\shortmid 2370
\newsymbol\shortparallel 2371
\newsymbol\smallsetminus 2272
\newsymbol\thicksim 2373
\newsymbol\thickapprox 2374
\newsymbol\approxeq 2375
\newsymbol\succapprox 2376
\newsymbol\precapprox 2377
\newsymbol\curvearrowleft 2378
\newsymbol\curvearrowright 2379
\newsymbol\digamma 207A
\newsymbol\varkappa 207B
\newsymbol\Bbbk 207C
\newsymbol\hslash 207D
\undefine\hbar
\newsymbol\hbar 207E
\newsymbol\backepsilon 237F
\catcode`\@=\csname pre amssym.tex at\endcsname


\centerline{\twelvebf Transformation de Fourier homog\`{e}ne}
\vskip 5mm
\centerline{by G\'{e}rard Laumon\footnote{${}^{\ast}$}{{\sevenrm Universit\'{e}
Paris-Sud et CNRS, UMR 8628, Math\'ematique, B\^{a}t.  425, F-91405
Orsay Cedex, France, Gerard.Laumon@math.u-psud.fr}}}
\vskip 20mm

Dans leur d\'{e}monstration de la correspondance de Drinfeld-Langlands
[1], Frenkel, Gaitsgory et Vilonen utilisent la transformation de
Fourier g\'{e}om\'{e}trique, ce qui les obligent soit \`{a} supposer
que le corps de base est de caract\'{e}ristique $0$ et \`{a}
travailler avec les ${\cal D}$-Modules, soit \`{a} supposer qu'il est
de caract\'{e}ristique $p>0$ et \`{a} travailler avec les faisceaux
$\ell$-adiques.  En fait, il n'utilisent cette transformation de
Fourier g\'{e}om\'{e}trique que pour des faisceaux homog\`{e}nes pour
lesquels on s'attend \`{a} avoir une transformation de Fourier sur
${\Bbb Z}$ qui prolonge la transformation de Radon g\'{e}om\'{e}trique
\'{e}tudi\'{e}e par Brylinski dans [2].

L'objet de cette note est de proposer une telle transformation de
Fourier.  Je remercie Bourbaki de m'avoir donn\'{e} l'occasion de
r\'{e}fl\'{e}chir \`{a} ce probl\`{e}me en pr\'{e}parant un
expos\'{e} \`{a} son s\'{e}minaire [3].  Je remercie aussi L. Illusie
pour son aide durant la r\'{e}daction de cette note.

\section{1}{La transformation de Fourier homog\`{e}ne}

1.1.  Soient $S$ un sch\'{e}ma et $\pi :V\rightarrow S$ un fibr\'{e}
vectoriel de rang constant $r$ sur $S$.  Le $S$-sch\'{e}ma en groupes
multiplicatif ${\Bbb G}_{{\rm m},S}$ agit par homoth\'{e}tie sur $V$
et on peut former le $S$-champ quotient
$$
V\,\smash{\mathop{\hbox to 6mm{\rightarrowfill}}\limits^{\scriptstyle
\rho}}\, {\cal V}=[V/{\Bbb G}_{{\rm m},S}]\,\smash{\mathop{\hbox to
6mm{\rightarrowfill}}\limits^{\scriptstyle\overline{\pi}}}\,S
$$
de $V$ par cette action.  Rappelons que, pour tout $S$-sch\'{e}ma $U$,
un objet de la cat\'{e}gorie ${\cal V}(U)$ est un couple $({\cal
L},x)$ form\'{e} d'un fibr\'{e} en droites ${\cal L}$ sur $U$ et d'un
morphisme de ${\cal O}_{U}$-Modules $x:{\cal L}\rightarrow {\cal
O}_{U}\otimes_{{\cal O}_{S}}V$ o\`{u} $V$ est ici consid\'{e}r\'{e}
comme un ${\cal O}_{S}$-Module localement libre de rang constant $r$.
Le $S$-champ ${\cal V}$ est alg\'{e}brique, de type fini et lisse
purement de dimension relative $r-1$, et le morphisme quotient $\rho$
est un ${\Bbb G}_{{\rm m},S}$-torseur repr\'{e}sentable.  Pour
$V={\Bbb A}_{S}^{1}$ le fibr\'{e} trivial de rang $1$, on note ${\cal
A}_{S}$ le champ ${\cal V}$ correspondant.

Soient $\pi^{\vee}:V^{\vee}\rightarrow S$ le fibr\'{e} vectoriel dual
de $\pi$ et
$$
\langle\,,\,\rangle :V^{\vee}\times_{S}V\rightarrow {\Bbb A}_{S}^{1}
$$
l'accouplement canonique.  Comme pr\'{e}c\'{e}demment, on forme le
$S$-champ quotient
$$
V^{\vee}\,\smash{\mathop{\hbox to 6mm{\rightarrowfill}}
\limits^{\scriptstyle\rho^{\vee}}}\, {\cal V}^{\vee}=[V^{\vee}/{\Bbb
G}_{{\rm m},S}]\,\smash{\mathop{\hbox to 6mm{\rightarrowfill}}
\limits^{\scriptstyle\overline{\pi}^{\vee}}}\,S
$$
et l'accouplement $\langle\,,\,\rangle$ induit un $S$-morphisme
$$
\mu :{\cal V}^{\vee}\times_{S}{\cal V}\rightarrow {\cal A}_{S}.
$$
Plus pr\'{e}cis\'{e}ment, si $U$ est un $S$-sch\'{e}ma, $\mu$ envoie
un objet $(({\cal M},y),({\cal L},x ))$ de $({\cal V}^{\vee}\times_{S}
{\cal V})(U)$ sur l'objet
$$
({\cal L}\otimes_{{\cal O}_{U}}{\cal M},\langle y ,x\rangle :{\cal L}
\otimes_{{\cal O}_{U}}{\cal M}\rightarrow {\cal O}_{U})
$$
de ${\cal A}_{S}(U)$.
\vskip 2mm

1.2.  On fixe un nombre premier $\ell$ inversible sur $S$ et une
cl\^{o}ture alg\'{e}brique $\overline{{\Bbb Q}}_{\ell}$ de ${\Bbb
Q}_{\ell}$.  On suppose que $S$ est de type fini sur un anneau
r\'{e}gulier $R$ de dimension $\leq 1$, de sorte que l'on dispose
d'un formalisme $\ell$-adique des six op\'{e}rations de Grothendieck
sur la cat\'{e}gorie des $S$-champs alg\'{e}briques de type fini (cf.
[4] et [5]).  Pour tout $S$-champ ${\cal X}$ de type fini, on dispose
donc des cat\'{e}gories d\'{e}riv\'{e}es $D_{{\rm c}}^{{\rm b}}({\cal
X}, \overline{{\Bbb Q}}_{\ell})\subset D_{{\rm c}}^{\pm}({\cal X},
\overline{{\Bbb Q}}_{\ell})$, avec:
\vskip 1mm

\itemitem{-} des foncteurs
$$
f_{\ast}:D_{{\rm c}}^{+}({\cal X},\overline{{\Bbb Q}}_{\ell})
\rightarrow D_{{\rm c}}^{+}({\cal Y},\overline{{\Bbb Q}}_{\ell}),
$$
$$
f_{!}:D_{{\rm c}}^{-}({\cal X},\overline{{\Bbb Q}}_{\ell})
\rightarrow D_{{\rm c}}^{-}({\cal Y},\overline{{\Bbb Q}}_{\ell}),
$$
$$
f^{\ast}:D_{{\rm c}}^{\pm ,{\rm b}}({\cal Y},\overline{{\Bbb
Q}}_{\ell})\rightarrow D_{{\rm c}}^{\pm ,{\rm b}}({\cal X},
\overline{{\Bbb Q}}_{\ell})
$$
et
$$
f^{!}:D_{{\rm c}}^{\pm ,{\rm b}}({\cal Y},\overline{{\Bbb Q}}_{\ell})
\rightarrow D_{{\rm c}}^{\pm ,{\rm b}}({\cal X},\overline{{\Bbb
Q}}_{\ell}),
$$
pour tout $S$-morphisme $f:{\cal X}\rightarrow {\cal Y}$ de $S$-champs
alg\'{e}briques de type fini, les foncteurs $f_{\ast}$ et $f_{!}$
respectant $D_{{\rm c}}^{{\rm b}}$ lorsque $f$ est repr\'{e}sentable,
\vskip 1mm

\itemitem{-} un foncteur {\og}{Hom interne}{\fg}
$$
\mathop{{\cal H}{\it om}}:D_{{\rm c}}^{{\rm b}}({\cal X},
\overline{{\Bbb Q}}_{\ell})^{{\rm opp}}\times D_{{\rm c}}^{\pm ,{\rm
b}}({\cal X},\overline{{\Bbb Q}}_{\ell})\rightarrow D_{{\rm c}}^{\pm
,{\rm b}} ({\cal X},\overline{{\Bbb Q}}_{\ell}),
$$
\vskip 1mm

\itemitem{-} un foncteur de dualit\'{e} involutif
$$
D_{{\cal X}}:D_{{\rm c}}^{\pm ,{\rm b}}({\cal X},\overline{{\Bbb
Q}}_{\ell})^{{\rm opp}}\rightarrow D_{{\rm c}}^{\mp ,{\rm b}}({\cal
X},\overline{{\Bbb Q}}_{\ell})
$$
d\'{e}fini par
$$
D_{{\cal X}}(K)=\mathop{{\cal H}{\it om}}(K,\varepsilon^{!}
\overline{{\Bbb Q}}_{\ell}[2d](d))
$$
o\`{u} $\varepsilon :{\cal X}\rightarrow\mathop{\rm Spec}(R)$ est le
morphisme structural et $d(=0\hbox{ ou }1)$ est la dimension de $R$,
\vskip 1mm

\itemitem{-} des isomorphismes canoniques de foncteurs $D_{{\cal
X}}\circ f_{\ast}\cong f_{!}\circ D_{{\cal X}}$ et $D_{{\cal X}}\circ
f^{\ast}\cong f^{!}\circ D_{{\cal X}}$,
\vskip 1mm

\itemitem{-} deux produits tensoriels internes \'{e}chang\'{e}s par la
dualit\'{e},
$$
\otimes :D_{{\rm c}}^{\pm ,{\rm b}}({\cal X},\overline{{\Bbb
Q}}_{\ell})\times D_{{\rm c}}^{\pm ,{\rm b}}({\cal X},
\overline{{\Bbb Q}}_{\ell})\rightarrow D_{{\rm c}}^{\pm ,{\rm b}}
({\cal X},\overline{{\Bbb Q}}_{\ell})
$$
et
$$\eqalign{
\widetilde{\otimes}:D_{{\rm c}}^{\pm ,{\rm b}}({\cal X},\overline{{\Bbb
Q}}_{\ell})\times D_{{\rm c}}^{\pm ,{\rm b}}({\cal X},
\overline{{\Bbb Q}}_{\ell})&\rightarrow D_{{\rm c}}^{\pm ,{\rm b}}
({\cal X},\overline{{\Bbb Q}}_{\ell}),\cr
(K_{1},K_{2})&\mapsto D_{{\cal X}}(D_{{\cal X}}(K_{1})\otimes
D_{{\cal X}}(K_{2})),\cr}
$$
\vskip 1mm

\itemitem{-} un isomorphisme canonique $\mathop{{\cal H}{\it
om}}(K,L)\cong D_{{\cal X}}(K)\widetilde{\otimes}L$ fonctoriel en
$(K,L)$ dans $D_{{\rm c}}^{{\rm b}}({\cal X},
\overline{{\Bbb Q}}_{\ell})^{{\rm opp}}\times D_{{\rm c}}^{\pm, {\rm
b}}({\cal X},\overline{{\Bbb Q}}_{\ell})$.
\vskip 1mm

On dispose aussi de la sous-cat\'{e}gorie strictement pleine
$\mathop{\rm Perv}({\cal X},\overline{{\Bbb Q}}_{\ell})\subset D_{{\rm
c}}^{{\rm b}}({\cal X},\overline{{\Bbb Q}}_{\ell})$ des
$\overline{{\Bbb Q}}_{\ell}$-faisceaux pervers pour la fonction de
dimension rectifi\'{e}e introduite par M. Artin (cf.  [6] et [7]).
Si $X$ est un sch\'{e}ma de type fini sur $S$, rappelons que cette
fonction est d\'{e}finie par
$$
\delta (x)=\mathop{\rm dim}(\overline{\{{\frak p}\}})+
\mathop{\rm deg.tr}(\kappa (x)/\kappa ({\frak p})),~
\forall x\in X,
$$
o\`{u} ${\frak p}\in \mathop{\rm Spec}(R)$ est l'image de $x$,
$\overline{\{{\frak p}\}}$ est l'adh\'{e}rence de $\{{\frak p}\}$ dans
$\mathop{\rm Spec}(R)$, $\kappa (x)$ est le corps r\'{e}siduel de $x$
et $\kappa ({\frak p})$ est le corps des fractions de $R/{\frak p}$,
et rappelons que $K\in\mathop{\rm ob}D_{{\rm c}}^{{\rm b}}
(X,\overline{{\Bbb Q}}_{\ell})$ est un $\overline{{\Bbb
Q}}_{\ell}$-faisceau pervers si et seulement si, pour tout point $x\in
X$, on a
$$
{\cal H}^{n}i_{x}^{\ast}K=(0),~\forall n>-\delta (x),\hbox{ et }
{\cal H}^{n}i_{x}^{!}K=(0),~\forall n<-\delta (x),
$$
o\`{u} $i_{x}:\{x\}\hookrightarrow X$ est l'inclusion.

La sous-cat\'{e}gorie $\mathop{\rm Perv}({\cal X},\overline{{\Bbb
Q}}_{\ell})\subset D_{{\rm c}}^{{\rm b}}({\cal X},\overline{{\Bbb
Q}}_{\ell})$ est stable par la dualit\'{e} $D_{{\cal X}}$ et elle est
noeth\'{e}rienne et artinienne (tous ses objets sont de longueur
finie).
\vskip 2mm

1.3.  L'inclusion naturelle ${\Bbb G}_{{\rm m},S}\hookrightarrow {\Bbb
A}_{S}^{1}$ passe au quotient en une immersion ouverte
$$
\beta :S=[{\Bbb G}_{{\rm m},S}/{\Bbb G}_{{\rm m},S}]\hookrightarrow
{\cal A}_{S}.
$$
On note
$$
\Psi =\beta_{\ast}\overline{{\Bbb Q}}_{\ell}\in D_{{\rm c}}^{{\rm
b}}({\cal A}_{S},\overline{{\Bbb Q}}_{\ell})
$$
l'image directe totale du faisceau constant $\overline{{\Bbb
Q}}_{\ell}$ sur $S$ par cette derni\`{e}re immersion ouverte.  Le
th\'{e}or\`{e}me de puret\'{e} cohomologique nous donne un triangle
distingu\'{e}
$$
\overline{{\Bbb Q}}_{\ell}\rightarrow\Psi\rightarrow
\alpha_{\ast}\overline{{\Bbb Q}}_{\ell}[-1](-1)\rightarrow
$$
dans $D_{{\rm c}}^{{\rm b}}({\cal A}_{S},\overline{{\Bbb Q}}_{\ell})$,
o\`{u} $\alpha :B({\Bbb G}_{{\rm m},S})=[S/{\Bbb G}_{{\rm
m},S}]\hookrightarrow {\cal A}_{S}$ est l'immersion ferm\'{e}e
compl\'{e}mentaire de $\beta$.

\thm LEMME 1.4
\enonce
{\rm (i)} Pour toute section $a:S\rightarrow {\Bbb A}_{S}^{1}$ de
${\Bbb A}_{S}^{1}$, on a
$$
h_{!}v_{\ast}\overline{{\Bbb Q}}_{\ell}=(0)
$$
o\`{u} $v:{\Bbb A}_{S}^{1}-a(S)\hookrightarrow {\Bbb A}_{S}^{1}$ est
l'ouvert compl\'{e}mentaire de l'image de la section et $h:{\Bbb
A}_{S}^{1}\rightarrow S$ est le morphisme structural.

\decale{\rm (ii)} On a un isomorphisme canonique
$$
\alpha^{\ast}\beta_{\ast}\overline{{\Bbb Q}}_{\ell}\buildrel\sim
\over\longrightarrow g_{!}\overline{{\Bbb Q}}_{\ell}[1]
$$
o\`{u} $g:S\rightarrow B({\Bbb G}_{{\rm m},S})$ est le ${\Bbb G}_{{\rm
m},S}$-torseur universel.
\endthm

\rem D\'{e}monstration
\endrem
Commen\c{c}ons par l'assertion (i).  Par le th\'{e}or\`{e}me de
changement de base propre, on peut supposer que $S$ est le spectre
d'un corps alg\'{e}briquement clos, et donc que $a$ est un point
ferm\'{e} de ${\Bbb A}^{1}$ (pour all\'{e}ger la r\'{e}daction, on
supprime l'indice $S$ des notations).  Il s'agit alors de
d\'{e}montrer que
$$
R\Gamma_{{\rm c}}({\Bbb A}^{1},v_{\ast}\overline{{\Bbb Q}}_{\ell})=
(0),
$$
ou ce qui revient au m\^{e}me par dualit\'{e}, que
$$
R\Gamma ({\Bbb A}^{1},v_{!}\overline{{\Bbb Q}}_{\ell})=
(0),
$$
Or, on a le triangle distingu\'{e}
$$
R\Gamma ({\Bbb A}^{1},v_{!}\overline{{\Bbb Q}}_{\ell})\rightarrow
R\Gamma ({\Bbb A}^{1},\overline{{\Bbb Q}}_{\ell})\rightarrow
\overline{{\Bbb Q}}_{\ell}\rightarrow
$$
o\`{u} la deuxi\`{e}me fl\`{e}che est la fl\`{e}che de restriction
\`{a} $\{a\}\subset {\Bbb A}^{1}$ et est donc un isomorphisme,
d'o\`{u} la conclusion.

D\'{e}montrons maintenant l'assertion (ii).  Pour cela consid\'{e}rons
le triangle distingu\'{e}
$$
\beta_{!}\overline{{\Bbb Q}}_{\ell}\rightarrow\beta_{\ast}
\overline{{\Bbb Q}}_{\ell}\rightarrow\alpha_{\ast}\alpha^{\ast}
\beta_{\ast}\overline{{\Bbb Q}}_{\ell}
$$
dans $D_{{\rm c}}^{{\rm b}}({\cal A}_{S},\overline{{\Bbb Q}}_{\ell})$ et
appliquons lui le foncteur $h_{!}$ pour le morphisme $h:{\cal A}_{S}=
[{\Bbb A}_{S}^{1}/{\Bbb G}_{{\rm m},S}]\rightarrow [S/{\Bbb G}_{{\rm
m},S}]=B({\Bbb G}_{{\rm m},S})$ qui est induit par le morphisme
structural ${\Bbb A}_{S}^{1}\rightarrow S$.  On obtient le triangle
distingu\'{e}
$$
g_{!}\overline{{\Bbb Q}}_{\ell}\rightarrow
h_{!}\beta_{\ast}\overline{{\Bbb Q}}_{\ell}\rightarrow
\alpha^{\ast}\beta_{\ast}\overline{{\Bbb Q}}_{\ell}
$$
dans $D_{{\rm c}}^{{\rm b}}(B({\Bbb G}_{{\rm m},S}),\overline{{\Bbb
Q}}_{\ell})$.  Ce dernier triangle d\'{e}g\'{e}n\`{e}re en un
isomorphisme
$$
\alpha^{\ast}\beta_{\ast}\overline{{\Bbb Q}}_{\ell}\buildrel\sim
\over\longrightarrow g_{!}\overline{{\Bbb Q}}_{\ell}[1]
$$
puisqu'on a $h_{!}\beta_{\ast}\overline{{\Bbb Q}}_{\ell}=(0)$
d'apr\`{e}s la partie (i) d\'{e}j\`{a} d\'{e}montr\'{e}e.
\hfill\hfill$\square$
\vskip 3mm

\thm D\'{E}FINITION 1.5
\enonce
La transformation de Fourier homog\`{e}ne est le foncteur
$$
\mathop{\rm Four}\nolimits_{{\cal V}/S}:D_{{\rm c}}^{{\rm b}}({\cal
V},\overline{{\Bbb Q}}_{\ell})\rightarrow D_{{\rm c}}^{-}({\cal
V}^{\vee},\overline{{\Bbb Q}}_{\ell})
$$
d\'{e}fini par
$$
\mathop{\rm Four}\nolimits_{{\cal V}/S}(K)=\mathop{\rm pr}
\nolimits_{!}^{\vee}(\mathop{\rm pr}\nolimits^{\ast}K
\otimes\mu^{\ast}\Psi )[r-1]
$$
o\`{u} $\mathop{\rm pr}\nolimits^{\vee}$ et $\mathop{\rm pr}$ sont les
projections canoniques de ${\cal V}^{\vee}\times_{S}{\cal V}$.
\endthm

Les immersions ouvertes $V^{\circ}\hookrightarrow V$ et $V^{\vee\circ}
\hookrightarrow V^{\vee}$ des compl\'{e}mentaires des sections nulles
des fibr\'{e}s vectoriels $V$ et $V^{\vee}$ passent au quotient en des
immersions ouvertes
$$
j:{\Bbb P}(V)={\cal V}^{\circ}=[V^{\circ}/{\Bbb G}_{{\rm m},S}]
\hookrightarrow {\cal V}\hbox{ et }j^{\vee}:{\Bbb P}(V^{\vee})=
{\cal V}^{\vee\circ}=[V^{\vee\circ}/{\Bbb G}_{{\rm m},S}]
\hookrightarrow {\cal V}^{\vee}
$$
o\`{u} ${\Bbb P}(V)$ et ${\Bbb P}(V^{\vee})$ sont les $S$-fibr\'{e}s
projectifs des droites de $V/S$ et $V^{\vee}/S$ respectivement.
Notons
$$
I:H\hookrightarrow {\Bbb P}(V^{\vee})\times_{S}{\Bbb P}(V)
$$
l'hypersurface d'incidence, quotient du ferm\'{e} de $V^{\vee\circ}
\times_{S}V^{\circ}$ d'\'{e}quation $\langle w,v\rangle =0$.  Notons
enfin $J:({\Bbb P}(V^{\vee})\times_{S}{\Bbb P}(V))- H\hookrightarrow
{\Bbb P}(V^{\vee})\times_{S}{\Bbb P}(V)$ l'immersion ouverte
compl\'{e}mentaire du ferm\'{e} $H$.  Comme $H$ est un diviseur lisse
sur $S$ dans un sch\'{e}ma lisse sur $S$, une nouvelle application du
th\'{e}or\`{e}me de puret\'{e} donne un triangle distingu\'{e}
$$
\overline{{\Bbb Q}}_{\ell}[2r-2]\rightarrow J_{\ast}\overline{{\Bbb
Q}}_{\ell}[2r-2]\rightarrow I_{\ast}\overline{{\Bbb
Q}}_{\ell}[2r-3](-1)\rightarrow
$$
dans $D_{{\rm c}}^{{\rm b}}({\Bbb P}(V^{\vee})\times_{S}{\Bbb P}(V),
\overline{{\Bbb Q}}_{\ell})$.

\thm PROPOSITION 1.6
\enonce
Le foncteur compos\'{e}
$$
(j^{\vee})^{\ast}\circ\mathop{\rm Four}\nolimits_{{\cal V}/S}\circ j_{!}:
D_{{\rm c}}^{{\rm b}}({\Bbb P}(V),\overline{{\Bbb Q}}_{\ell})\rightarrow
D_{{\rm c}}^{-}({\Bbb P}(V^{\vee}),\overline{{\Bbb Q}}_{\ell})
$$
n'est autre que le foncteur
$$
K\mapsto\mathop{\rm pr}\nolimits_{!}^{\vee}(\mathop{\rm pr}
\nolimits^{\ast}K\otimes J_{\ast}\overline{{\Bbb Q}}_{\ell})[r-1]
$$
o\`{u} $\mathop{\rm pr}^{\vee}$ et $\mathop{\rm pr}$ sont maintenant
les deux projections canoniques de ${\Bbb P}(V^{\vee})\times_{S}{\Bbb
P}(V)$, et on a un d\'{e}vissage fonctoriel en $K\in\mathop{\rm ob}D_{{\rm
c}}^{{\rm b}}({\Bbb P}(V),\overline{{\Bbb Q}}_{\ell})$,
$$
(\overline{\pi}^{\vee\circ})^{\ast}(\overline{\pi}^{\circ})_{!}K[r-1]
\rightarrow (j^{\vee})^{\ast} \mathop{\rm Four}\nolimits_{{\cal V}/S}
(j_{!}K) \rightarrow\mathop{\rm Rad}\nolimits_{{\Bbb P}(V)/S}
(K)(-1)\rightarrow
$$
o\`{u} $\overline{\pi}^{\vee\circ}:{\Bbb P}(V^{\vee})\rightarrow S$ et
$\overline{\pi}^{\circ}:{\Bbb P}(V)\rightarrow S$ sont les morphismes
structuraux et o\`{u}
$$
\mathop{\rm Rad}\nolimits_{{\Bbb P}(V)/S}:D_{{\rm c}}^{{\rm b}}({\Bbb
P}(V),\overline{{\Bbb Q}}_{\ell})\rightarrow D_{{\rm c}}^{{\rm
b}}({\Bbb P}(V^{\vee}),\overline{{\Bbb Q}}_{\ell})
$$
est la transformation de Radon g\'{e}om\'{e}trique {\rm (}cf.  {\rm
[2])} d\'{e}finie par
$$
\mathop{\rm Rad}\nolimits_{{\Bbb P}(V)/S}(K)=q_{!}q^{\ast}K[r-2]
$$
o\`{u} $q^{\vee}:H\rightarrow {\Bbb P}(V^{\vee})$ et $q:H\rightarrow {\Bbb
P}(V)$ sont les restrictions \`{a} $H$ des deux projections canoniques.
\endthm

\rem D\'{e}monstration
\endrem
Comme la restriction du morphisme $\mu : {\cal V}^{\vee}\times_{S}{\cal
V}\rightarrow {\cal A}$ \`{a} l'ouvert ${\Bbb P}(V^{\vee})\times_{S}{\Bbb
P}(V)$ est lisse, on a
$$
(j^{\vee}\times_{S}j)^{\ast}\mu^{\ast}\Psi\cong
J_{\ast}\overline{{\Bbb Q}}_{\ell}
$$
par le th\'{e}or\`{e}me de changement de base par un morphisme lisse,
d'o\`{u} la premi\`{e}re assertion.

Compte tenu de cette premi\`{e}re assertion, on d\'{e}duit du
triangle distingu\'{e} pr\'{e}c\'{e}dant l'\'{e}nonc\'{e} de la
proposition un triangle distingu\'{e}
$$
\mathop{\rm pr}\nolimits_{!}^{\vee}\mathop{\rm pr}
\nolimits^{\ast}K[r-1]\rightarrow (j^{\vee})^{\ast}\mathop{\rm Four}
\nolimits_{{\cal V}/S}(j_{!}K)\rightarrow \mathop{\rm pr}
\nolimits_{!}^{\vee}(\mathop{\rm pr}\nolimits^{\ast}K\otimes
I_{\ast}\overline{{\Bbb Q}}_{\ell})[r-2](-1)\rightarrow ,
$$
d'o\`{u} la seconde assertion puisque $\mathop{\rm pr}
\nolimits_{!}^{\vee}\mathop{\rm pr} \nolimits^{\ast}\cong
(\overline{\pi}^{\vee\circ})^{\ast} (\overline{\pi}^{\circ})_{!}$
d'apr\`{e}s le th\'{e}or\`{e}me de changement de base propre.
\hfill\hfill$\square$
\vskip 3mm

Soient $W\rightarrow S$ un deuxi\`{e}me fibr\'{e} vectoriel de rang
constant $s$ et $f:V\rightarrow W$ un morphisme de $S$-fibr\'{e}s
vectoriels de transpos\'{e} ${}^{{\rm t}}\!f:W^{\vee}\rightarrow V^{\vee}$.
Par passage au quotient on a des $S$-morphismes, not\'{e}s encore $f$ et
${}^{{\rm t}}\!f$,
$$
f:{\cal V}\rightarrow {\cal W}\hbox{ et }{}^{{\rm t}}\!f:{\cal
W}^{\vee}\rightarrow {\cal V}^{\vee}.
$$
et on v\'{e}rifie, comme on le fait pour la transformation de
Fourier-Deligne [8], le lemme suivant.

\thm LEMME 1.7
\enonce
On a un isomorphisme canonique de foncteurs
$$
\mathop{\rm Four}\nolimits_{{\cal W}/S}\circ f_{!}\cong
({}^{{\rm t}}\!f)^{\ast}\circ\mathop{\rm Four}\nolimits_{{\cal V}/S}[s-r].
\eqno{\square}
$$
\endthm

En particulier, si $i:B({\Bbb G}_{{\rm m},S})\hookrightarrow {\cal V}$
et $i^{\vee}:B({\Bbb G}_{{\rm m},S})\hookrightarrow {\cal V}^{\vee}$
sont les immersions ferm\'{e}es compl\'{e}mentaires des immersions
ouvertes $j$ et $j^{\vee}$ et si on note encore $\pi :{\cal
V}\rightarrow B({\Bbb G}_{{\rm m},S})$ et $\pi^{\vee}:{\cal
V}^{\vee}\rightarrow B({\Bbb G}_{{\rm m},S})$ les morphismes
repr\'{e}sentables induits par les projections canoniques $\pi
:V\rightarrow S$ et $\pi^{\vee}:V^{\vee}\rightarrow S$, on a
$$
\mathop{\rm Four}\nolimits_{{\cal V}/S}\circ i_{!}\cong
(\pi^{\vee})^{\ast}\circ\mathop{\rm Four}\nolimits_{B({\Bbb G}_{{\rm
m},S})/S}[r]
$$
et
$$
\mathop{\rm Four}\nolimits_{B({\Bbb G}_{{\rm m},S})/S}
\circ\pi_{!}\cong (i^{\vee})^{\ast}\circ\mathop{\rm Four}
\nolimits_{{\cal V}/S}[-r].
$$

\thm PROPOSITION 1.8
\enonce
La transformation de Fourier homog\`{e}ne
$$
\mathop{\rm Four}\nolimits_{B({\Bbb G}_{{\rm m},S})/S}:D_{{\rm
c}}^{{\rm b}}(B({\Bbb G}_{{\rm m},S}),\overline{{\Bbb Q}}_{\ell})
\rightarrow D_{{\rm c}}^{-}(B({\Bbb G}_{{\rm m},S}),\overline{{\Bbb
Q}}_{\ell})
$$
est le foncteur d'inclusion de $D_{{\rm c}}^{{\rm b}}(B({\Bbb G}_{{\rm
m},S}),\overline{{\Bbb Q}}_{\ell})$ dans $D_{{\rm c}}^{-} (B({\Bbb
G}_{{\rm m},S}),\overline{{\Bbb Q}}_{\ell})$.
\endthm

\rem D\'{e}monstration
\endrem
Par d\'{e}finition, on a
$$
\mathop{\rm Four}\nolimits_{B({\Bbb G}_{{\rm m},S})/S}(K)=\mathop{\rm
pr}\nolimits_{!}^{\vee}(\mathop{\rm pr}\nolimits^{\ast}K\otimes
\nu^{\ast}\alpha^{\ast}\Psi )[-1]
$$
o\`{u} $\mathop{\rm pr}\nolimits^{\vee}$ et $\mathop{\rm pr}$ sont les
projections canoniques de $B({\Bbb G}_{{\rm m},S})\times_{S}B({\Bbb
G}_{{\rm m},S})$ et o\`{u} $\nu :B({\Bbb G}_{{\rm m},S})\times_{S}
B({\Bbb G}_{{\rm m},S})\rightarrow B({\Bbb G}_{{\rm m},S})$ envoie
tout couple de fibr\'{e}s en droites $({\cal M},{\cal L})$ sur le
fibr\'{e} en droites ${\cal L}\otimes {\cal M}^{\otimes -1}$.  De
plus, on a un isomorphisme
$$
\alpha^{\ast}\Psi\buildrel\sim\over\longrightarrow g_{!}
\overline{{\Bbb Q}}_{\ell}[1]
$$
d'apr\`{e}s 1.4(ii).  Enfin, on a le carr\'{e} cart\'{e}sien
$$\diagram{
B({\Bbb G}_{{\rm m},S})&\kern -10mm\smash{\mathop{\hbox to
22mm{\rightarrowfill}}\limits^{\scriptstyle\varepsilon}}\kern
-6mm&S\cr
\llap{$\scriptstyle\Delta$}\left\downarrow
\vbox to 4mm{}\right.\rlap{}&\square&\llap{}\left\downarrow
\vbox to 4mm{}\right.\rlap{$\scriptstyle g$}\cr
B({\Bbb G}_{{\rm m},S})\times_{S}B({\Bbb G}_{{\rm m},S})&\kern
-1mm\smash{\mathop{\hbox to 8mm{\rightarrowfill}}
\limits_{\scriptstyle\nu}}\kern -1mm&B({\Bbb G}_{{\rm m},S})\cr}
$$
o\`{u} $\Delta$ est le morphisme morphisme diagonal
et $\varepsilon$ est le morphisme structural.

Appliquant le th\'{e}or\`{e}me de changement de base propre \`{a} ce
carr\'{e} et la formule des projections, on obtient que
$$
\mathop{\rm pr}\nolimits_{!}^{\vee}(\mathop{\rm pr}
\nolimits^{\ast}K\otimes\nu^{\ast}g_{!}\overline{{\Bbb Q}}_{\ell})=
\mathop{\rm pr}\nolimits_{!}^{\vee}(\mathop{\rm pr}
\nolimits^{\ast}K\otimes\Delta_{!}\overline{{\Bbb Q}}_{\ell})=K,
$$
d'o\`{u} la proposition.
\hfill\hfill$\square$
\vskip 3mm

\thm COROLLAIRE 1.9
\enonce
L'image essentielle de la transformation de Fourier homog\`{e}ne
$\mathop{\rm Four}\nolimits_{{\cal V}/S}$ est contenue dans la
sous-cat\'{e}gorie pleine $D_{{\rm c}}^{{\rm b}}({\cal
V}^{\vee}/S,\overline{{\Bbb Q}}_{\ell})$ de $D_{{\rm c}}^{-}({\cal
V}^{\vee}/S,\overline{{\Bbb Q}}_{\ell})$.
\endthm

\rem D\'{e}monstration
\endrem
Tout $K\in\mathop{\rm ob} D_{{\rm c}}^{{\rm b}}({\cal
V}/S,\overline{{\Bbb Q}}_{\ell})$ admet le d\'{e}vissage
$$
j_{!}j^{\ast}K\rightarrow K\rightarrow i_{\ast}i^{\ast}K\rightarrow ,
$$
de sorte que $\mathop{\rm Four}\nolimits_{{\cal V}/S}(K)$ admet le
d\'{e}vissage
$$
\mathop{\rm Four}\nolimits_{{\cal V}/S}(j_{!}j^{\ast}K)\rightarrow
\mathop{\rm Four}\nolimits_{{\cal V}/S}(K)\rightarrow\mathop{\rm
Four}\nolimits_{{\cal V}/S}(i_{!}i^{\ast}K)\rightarrow
$$
et que $\mathop{\rm Four}\nolimits_{{\cal V}/S}(j_{!}j^{\ast}K)$ admet
\`{a} son tour le d\'{e}vissage
$$
j_{!}^{\vee}(j^{\vee})^{\ast}\mathop{\rm Four}\nolimits_{{\cal V}/S}
(j_{!}j^{\ast}K)\rightarrow\mathop{\rm Four}\nolimits_{{\cal V}/S}
(j_{!}j^{\ast}K)\rightarrow i_{\ast}^{\vee}(i^{\vee})^{\ast}
\mathop{\rm Four}\nolimits_{{\cal V}/S}(j_{!}j^{\ast}K)\rightarrow .
$$
D'apr\`{e}s la proposition 1.6, $(j^{\vee})^{\ast}\mathop{\rm
Four}\nolimits_{{\cal V}/S}(j_{!}j^{\ast}K)$ est \`{a} cohomologie
born\'{e}e et, d'apr\`{e}s le lemme 1.7 et la proposition 1.8, il en
est de m\^{e}me de
$$
\mathop{\rm Four}\nolimits_{{\cal V}/S}(i_{!}i^{\ast}K)=
(\pi^{\vee})^{\ast}\mathop{\rm Four}\nolimits_{B({\Bbb G}_{{\rm
m},S})/S} (i^{\ast}K)[r]=(\pi^{\vee})^{\ast}i^{\ast}K[r]
$$
et de
$$
(i^{\vee})^{\ast}\mathop{\rm Four}\nolimits_{{\cal V}/S}(j_{!}j^{\ast}K)=
\mathop{\rm Four}\nolimits_{B({\Bbb G}_{{\rm m},S})/S}
(\pi_{!}j_{!}j^{\ast}K) =\pi_{!}^{\circ}j^{\ast}K[r]
$$
o\`{u} $\pi^{\circ}=\pi\circ j:{\Bbb P}(V)\rightarrow B({\Bbb G}_{{\rm
m},S})$, d'o\`{u} le corollaire.
\hfill\hfill$\square$

\section{2}{Relation entre les transformations de Fourier homog\`{e}ne
et de Fourier-Deligne}

2.1.  Supposons dans cette section que $S$ est de caract\'{e}ristique
$p>0$ et fixons un caract\`{e}re additif non trivial $\psi :{\Bbb
F}_{p}\rightarrow\overline{{\Bbb Q}}_{\ell}^{\times}$.

La transformation de Fourier-Deligne associ\'{e}e \`{a} $V/S$ et
\`{a} $\psi$ (cf. [8] et [9]) est le foncteur
$$
\mathop{\rm Four}\nolimits_{V/S,\psi}:D_{{\rm c}}^{{\rm
b}}(V,\overline{{\Bbb Q}}_{\ell})\rightarrow D_{{\rm c}}^{{\rm
b}}(V^{\vee},\overline{{\Bbb Q}}_{\ell})
$$
d\'{e}fini par
$$
\mathop{\rm Four}\nolimits_{V/S,\psi}(K)=\mathop{\rm pr}
\nolimits_{!}^{\vee}(\mathop{\rm pr}\nolimits^{\ast}K
\otimes\langle\,,\,\rangle ^{\ast}{\cal L}_{\psi})[r]
$$
o\`{u} $\mathop{\rm pr}\nolimits^{\vee}$ et $\mathop{\rm pr}$ sont les
projections canoniques de $V^{\vee}\times_{S}V$, o\`{u}
$$
\langle\,,\,\rangle : V^{\vee}\times_{S}V\rightarrow {\Bbb
A}_{S}^{1}
$$
est l'accouplement canonique et o\`{u} ${\cal L}_{\psi}$ est le
$\overline{{\Bbb Q}}_{\ell}$-faisceau lisse de rang $1$
d'Artin-Schreier sur ${\Bbb A}_{S}^{1}$ associ\'{e} au
caract\`{e}re $\psi$.

\thm TH\'{E}OR\`{E}ME 2.2
\enonce
Les foncteurs compos\'{e}s
$$
(\rho^{\vee})^{\ast}\circ\mathop{\rm Four}\nolimits_{{\cal V}/S}\hbox{
et }\mathop{\rm Four}\nolimits_{V/S,\psi}\circ\rho^{\ast}:D_{{\rm
c}}^{{\rm b}}({\cal V},\overline{{\Bbb Q}}_{\ell})\rightarrow D_{{\rm
c}}^{{\rm b}}(V^{\vee},\overline{{\Bbb Q}}_{\ell})
$$
sont canoniquement isomorphes.
\endthm

\rem D\'{e}monstration
\endrem
Faisons agir ${\Bbb G}_{{\rm m},S}$ sur $V^{\vee}\times_{S}V$ de
mani\`{e}re anti-diagonale ($t\cdot (w,v)=(t^{-1}w,tv)$).  On
a le diagramme \`{a} carr\'{e}s cart\'{e}siens
$$\diagram{
V^{\vee}&\kern -1mm\smash{\mathop{\hbox to 15mm{\leftarrowfill}}
\limits^{\scriptstyle\mathop{\rm pr}\nolimits^{\vee}}}\kern -8mm&
V^{\vee}\times_{S}V&\kern -8mm\smash{\mathop{\hbox to 15mm
{\rightarrowfill}}\limits^{\scriptstyle\mathop{\rm pr}}}\kern
-1mm&V\cr
\llap{$\scriptstyle\rho^{\vee}$}\left\downarrow\vbox to
4mm{}\right.\rlap{}&\square&\llap{$\scriptstyle {\rm R}$}\left\downarrow
\vbox to 4mm{}\right.\rlap{}&\square&\llap{}\left\downarrow\vbox to
4mm{}\right.\rlap{$\scriptstyle\rho$}\cr
{\cal V}^{\vee}&\kern -1mm\smash{\mathop{\hbox to 8mm{\leftarrowfill}}
\limits_{\scriptstyle\mathop{\overline{\rm pr}}\nolimits^{\vee}}}\kern -1mm
&[(V^{\vee}\times_{S}V)/{\Bbb G}_{{\rm m},S}] &\kern
-1mm\smash{\mathop{\hbox to 8mm{\rightarrowfill}}
\limits_{\scriptstyle\mathop{\overline{\rm pr}}}}\kern -1mm&{\cal V}\cr}
$$
o\`{u} $R$ est le morphisme quotient, et $\langle\,,\,\rangle :
V^{\vee}\times_{S}V\rightarrow {\Bbb A}_{S}^{1}$ se factorise en
$$
V^{\vee}\times_{S}V\,\smash{\mathop{\hbox to 6mm{\rightarrowfill}}
\limits^{\scriptstyle\mathop{\rm R}}}\,
[(V^{\vee}\times_{S}V)/{\Bbb G}_{{\rm m},S}]\,\smash{\mathop{\hbox to
8mm{\rightarrowfill}}\limits^{\scriptstyle\overline{\langle\,,\,\rangle}}}
\,{\Bbb A}_{S}^{1}
$$

Une application du th\'{e}or\`{e}me de changement de base pour un
morphisme repr\'{e}sentable propre montre que, pour tout
$K\in\mathop{\rm ob}D_{{\rm c}}^{{\rm b}}({\cal V},\overline{{\Bbb
Q}}_{\ell})$
$$
\mathop{\rm Four}\nolimits_{V/S,\psi}(\rho^{\ast}K)=(\rho^{\vee})^{\ast}
\mathop{\overline{\rm pr}}\nolimits_{!}^{\vee}(\mathop{\overline{\rm pr}}
\nolimits^{\ast}K\otimes\overline{\langle\,,\,\rangle}^{\ast}{\cal
L}_{\psi})[r].
$$
En termes imag\'{e}s, le transform\'{e} de Fourier-Deligne d'un
complexe ${\Bbb G}_{{\rm m},S}$-\'{e}quivariant est aussi ${\Bbb
G}_{{\rm m},S}$-\'{e}quivariant.

Maintenant, ${\cal V}^{\vee}\times_{S}{\cal V}$ est le quotient de
$[(V^{\vee}\times_{S}V)/{\Bbb G}_{{\rm m},S}]$ par l'action de ${\Bbb
G}_{{\rm m},S}$ induite par celle par homoth\'{e}tie sur le facteur
$V$ de $V^{\vee}\times_{S}V$. Si on note
$$
F:[(V^{\vee}\times_{S}V)/{\Bbb G}_{{\rm m},S}]\rightarrow {\cal
V}^{\vee}\times_{S}{\cal V}
$$
le morphisme quotient, les morphismes $\mathop{\overline{\rm pr}}$ et
$\mathop{\overline{\rm pr}}\nolimits^{\vee}$ introduits ci-dessus se
factorisent par $F$ en les projections canoniques $\mathop{\rm pr}$ et
$\mathop{\rm pr}\nolimits^{\vee}$ de ${\cal V}^{\vee}\times_{S} {\cal
V}$.  Par suite, une application de la formule des projections nous
donne un isomorphisme canonique
$$
\mathop{\overline{\rm pr}}\nolimits_{!}^{\vee}(\mathop{\overline{\rm
pr}}\nolimits^{\ast}K\otimes\overline{\langle\,,\,\rangle}^{\ast}{\cal
L}_{\psi})[r]\cong\mathop{\rm pr}\nolimits_{!}^{\vee}(\mathop{\rm pr}
\nolimits^{\ast}K\otimes F_{!}\overline{\langle\,,\,\rangle}^{\ast}{\cal
L}_{\psi})[r].
$$

Calculons enfin $F_{!}\overline{\langle\,,\,\rangle}^{\ast}{\cal
L}_{\psi}$.  On a le carr\'{e} cart\'{e}sien
$$\diagram{
[(V^{\vee}\times_{S}V)/{\Bbb G}_{{\rm m},S}]&\kern
-1mm\smash{\mathop{\hbox to 8mm{\rightarrowfill}}
\limits^{\scriptstyle\overline{\mu}}}\kern -2mm&{\Bbb A}_{S}^{1}\cr
\llap{$\scriptstyle F$}\left\downarrow\vbox to
4mm{}\right.\rlap{}&\square&\llap{}\left\downarrow\vbox to
4mm{}\right.\rlap{$\scriptstyle f$}\cr
{\cal V}^{\vee}\times_{S}{\cal V}&\kern -8mm\smash{\mathop{\hbox to
15mm{\rightarrowfill}}\limits_{\scriptstyle\mu}}\kern -2mm&{\cal
A}_{S}\cr}
$$
o\`{u} $f:{\Bbb A}_{S}^{1}\rightarrow [{\Bbb A}_{S}^{1}/{\Bbb G}_{{\rm
m},S}]={\cal A}_{S}$ est le morphisme quotient.  En appliquant de
nouveau le th\'{e}or\`{e}me de changement de base propre, on obtient
donc un isomorphisme canonique
$$
F_{!}\overline{\langle\,,\,\rangle}^{\ast}{\cal L}_{\psi}\cong
\mu^{\ast}f_{!}{\cal L}_{\psi}.
$$
On conclut par le lemme 2.3 ci-dessous.
\hfill\hfill$\square$
\vskip 3mm

\thm LEMME 2.3
\enonce
Le complexe $\beta^{\ast}f_{!}{\cal L}_{\psi}$ est canoniquement
isomorphe au faisceau constant $\overline{{\Bbb Q}}_{\ell}$ sur $S$,
plac\'{e} en degr\'{e} $1$, et le morphisme d'adjonction $f_{!}{\cal
L}_{\psi}\rightarrow \beta_{\ast}\beta^{\ast}f_{!}{\cal L}_{\psi}$ est
un isomorphisme.  Par suite, le complexe $f_{!}{\cal L}_{\psi}[1]$ est
canoniquement isomorphe au complexe $\Psi =\beta_{\ast}\overline{{\Bbb
Q}}_{\ell}$ d\'{e}fini en $1.3$.
\endthm

\rem D\'{e}monstration
\endrem
Si $h:{\Bbb A}_{S}^{1}\rightarrow S$ est le morphisme structural, il
est bien connu que $h_{!}{\cal L}_{\psi}=(0)$.  Par suite, 
$\beta^{\ast}f_{!}{\cal
L}_{\psi}$, qui est l'image directe par le morphisme structural
${\Bbb G}_{{\rm m},S}\rightarrow S$ de la restriction \`{a} ${\Bbb
G}_{{\rm m},S}\subset {\Bbb A}_{S}^{1}$ de ${\cal L}_{\psi}$, est
canoniquement isomorphe \`{a} $\overline{{\Bbb Q}}_{\ell}[-1]$.

Le c\^{o}ne de la fl\`{e}che d'adjonction $f_{!}{\cal L}_{\psi}
\rightarrow \beta_{\ast}\beta^{\ast}f_{!}{\cal L}_{\psi}$ est le
complexe $\alpha^{!}f_{!}{\cal L}_{\psi}[1]$ et il nous reste \`{a}
d\'{e}montrer que $\alpha^{!}f_{!}{\cal L}_{\psi}=(0)$, ou ce qui
revient au m\^{e}me par dualit\'{e}, que $\alpha^{\ast}f_{\ast}{\cal
L}_{\psi^{-1}}=(0)$.  Pour cela, on peut se restreindre \`{a} ${\Bbb
A}_{S}^{1}$ par le morphisme lisse $f:{\Bbb A}_{S}^{1}\rightarrow
{\cal A}_{S}$.  Comme le carr\'{e}
$$\diagram{
{\Bbb A}_{S}^{1}\times_{S}{\Bbb G}_{{\rm m},S}&\kern
-1mm\smash{\mathop{\hbox to 12mm{\rightarrowfill}}
\limits^{\scriptstyle (x,t)\mapsto xt}}\kern -1mm&{\Bbb A}_{S}^{1}\cr
\llap{$\scriptstyle \mathop{\rm pr}$}\left\downarrow
\vbox to 4mm{}\right.\rlap{}&\square&\llap{}\left\downarrow
\vbox to 4mm{}\right.\rlap{$\scriptstyle f$}\cr
{\Bbb A}_{S}^{1}&\kern -8mm\smash{\mathop{\hbox to
19mm{\rightarrowfill}} \limits_{\scriptstyle f}}\kern -1mm&{\cal
A}_{S}\cr}
$$
est cart\'{e}sien, une application du th\'{e}or\`{e}me de changement
de base lisse nous ram\`{e}ne donc \`{a} v\'{e}rifier que
$0^{\ast}\mathop{\rm pr}\nolimits_{\ast}{\cal L}_{\psi^{-1}}(xt)=(0)$,
o\`{u} $\mathop{\rm pr}$ est la projection canonique sur le premier
facteur de ${\Bbb A}_{S}^{1}\times_{S}{\Bbb G}_{{\rm m},S}$ et
$0:S\rightarrow {\Bbb A}_{S}^{1}$ est la section nulle.

Compactifions $\mathop{\rm pr}$ en
$$\diagram{
{\Bbb A}_{S}^{1}\times_{S}{\Bbb G}_{{\rm m},S}&\kern -2mm
\smash{\mathop{\lhook\joinrel\mathrel{\hbox to 6mm{\rightarrowfill}}
}\limits^{\scriptstyle }}\kern -2mm&{\Bbb A}_{S}^{1}\times_{S}{\Bbb 
P}_{S}^{1}\cr
\llap{$\scriptstyle \mathop{\rm pr}$}\left\downarrow
\vbox to 5mm{}\right.\rlap{}&&\cr
{\Bbb A}_{S}^{1}&&\cr
\arrow(28,16)\dir(-3,-2)\length{13}
\put (27,10){\scriptstyle\mathop{\overline{\rm pr}}}}
$$
le ferm\'{e} compl\'{e}mentaire de ${\Bbb A}_{S}^{1}\times_{S}{\Bbb
G}_{{\rm m},S}$ dans ${\Bbb A}_{S}^{1}\times_{S}{\Bbb P}_{S}^{1}$
\'{e}tant la r\'{e}union disjointe des images des deux sections
\'{e}videntes ${\Bbb A}_{S}^{1} \rightarrow {\Bbb
A}_{S}^{1}\times_{S}{\Bbb P}_{S}^{1},~x\mapsto (x,0)\hbox{ et
}x\mapsto (x,\infty)$.

Rappelons que le syst\`{e}me local ${\cal L}_{\psi^{-1}}(xt)$ sur
${\Bbb A}_{S}^{1}\times_{S}{\Bbb G}_{{\rm m},S}\subset {\Bbb
A}_{S}^{1}\times_{S}{\Bbb P}_{S}^{1}$ se prolonge en un syst\`{e}me
local, encore not\'{e} ${\cal L}_{\psi^{-1}}(xt)$, sur ${\Bbb
A}_{S}^{1}\times_{S}{\Bbb A}_{S}^{1}\subset {\Bbb A}_{S}^{1}
\times_{S}{\Bbb P}_{S}^{1}$, et que les prolongements par z\'{e}ro et
par image directe totale de ce dernier syst\`{e}me local \`{a} ${\Bbb
A}_{S}^{1} \times_{S}{\Bbb P}_{S}^{1}$ tout entier co\"{\i}ncident
(cf.  la d\'{e}monstration de (2.4.1) dans [9]).

Si on note $\overline{{\cal L}}_{\psi^{-1}}(xt)$
ces deux m\^{e}mes prolongements, on a un triangle distingu\'{e}
$$
\mathop{\overline{\rm pr}}\nolimits_{\ast}\overline{{\cal
L}}_{\psi^{-1}}(xt)\rightarrow \mathop{\rm pr}\nolimits_{\ast}{\cal
L}_{\psi^{-1}}(xt)\rightarrow \overline{{\Bbb Q}}_{\ell}[-1](-1)
\rightarrow
$$
dans $D_{{\rm c}}^{{\rm b}}({\Bbb A}_{S}^{1},\overline{{\Bbb
Q}}_{\ell})$, et donc un triangle distingu\'{e}
$$
0^{\ast}\mathop{\overline{\rm pr}}\nolimits_{\ast}\overline{{\cal
L}}_{\psi^{-1}}(xt)\rightarrow 0^{\ast}\mathop{\rm pr}\nolimits_{\ast}
{\cal L}_{\psi^{-1}}(xt)\rightarrow \overline{{\Bbb
Q}}_{\ell}[-1](-1)\rightarrow
$$
dans $D_{{\rm c}}^{{\rm b}}(S,\overline{{\Bbb Q}}_{\ell})$,
et il ne reste plus qu'\`{a} v\'{e}rifier que le fl\`{e}che de
bord
$$
\overline{{\Bbb Q}}_{\ell}[-1](-1)\rightarrow 0^{\ast}
\mathop{\overline{\rm pr}}\nolimits_{\ast}\overline{{\cal
L}}_{\psi^{-1}}(xt)[1]
$$
est un isomorphisme.  Mais, par le th\'{e}or\`{e}me de changement de
base propre, cette fl\`{e}che n'est autre, \`{a} un d\'{e}calage
pr\`{e}s, que l'inverse de l'isomorphisme trace
$$
h_{!}\overline{{\Bbb Q}}_{\ell}\buildrel\sim\over\longrightarrow
\overline{{\Bbb Q}}_{\ell}[-2](-1),
$$
d'o\`{u} la conclusion.
\hfill\hfill$\square$
\vskip 3mm

2.4. Il r\'{e}sulte du th\'{e}or\`{e}me $2.1$ que, si $S$ est de
caract\'{e}ristique $p>0$, la transformation de Fourier homog\`{e}ne
est involutive, commute \`{a} la dualit\'{e} et induit une
\'{e}quivalence de la cat\'{e}gorie des $\overline{{\Bbb
Q}}_{\ell}$-faisceaux pervers sur ${\cal V}$ sur celle des
$\overline{{\Bbb Q}}_{\ell}$-faisceaux pervers sur ${\cal V}^{\vee}$.
Dans les deux sections qui suivent, nous allons donner une
d\'{e}monstration directe de ces propri\'{e}t\'{e}s et donc voir
qu'elles valent aussi pour $S$ arbitraire.

\section{3}{Involutivit\'{e} de la transformation de Fourier
homog\`{e}ne}

\thm TH\'{E}OR\`{E}ME 3.1
\enonce
La transformation de Fourier homog\`{e}ne est involutive au sens
o\`{u} le foncteur compos\'{e}
$$
\mathop{\rm Four}\nolimits_{{\cal V}^{\vee}/S}\circ\mathop{\rm
Four}\nolimits_{{\cal V}/S}:D_{{\rm c}}^{{\rm b}}({\cal V},\overline{{\Bbb
Q}}_{\ell})\rightarrow D_{{\rm c}}^{{\rm b}}({\cal V},\overline{{\Bbb
Q}}_{\ell})
$$
est canoniquement isomorphe au foncteur $K\mapsto K(-r)$.
\endthm

Avant de d\'{e}montrer le th\'{e}or\`{e}me, \'{e}tablissons quelques
propri\'{e}t\'{e}s du complexe $\Psi =\beta_{\ast}\overline{{\Bbb
Q}}_{\ell}\in D_{{\rm c}}^{{\rm b}}({\cal A}_{S},\overline{{\Bbb
Q}}_{\ell})$.

\thm LEMME 3.2
\enonce
Soient $f:{\Bbb A}_{S}^{1}\rightarrow {\cal A}_{S}$ le morphisme quotient et
$v:{\Bbb A}_{S}^{1}-a(S)\hookrightarrow {\Bbb A}_{S}^{1}$ l'ouvert 
compl\'{e}mentaire de
l'image d'une section $a:S\rightarrow {\Bbb A}_{S}^{1}$ {\rm partout 
non nulle}.
Le complexe $\beta^{\ast}f_{!}v_{\ast}\overline{{\Bbb Q}}_{\ell}$ est
canoniquement isomorphe au faisceau constant $\overline{{\Bbb
Q}}_{\ell}$ sur $S$, plac\'{e} en degr\'{e} $1$, et le morphisme
d'adjonction $f_{!}v_{\ast}\overline{{\Bbb Q}}_{\ell}\rightarrow
\beta_{\ast}\beta^{\ast}f_{!}v_{\ast}\overline{{\Bbb Q}}_{\ell}$ est
un isomorphisme.  Par suite, le complexe $f_{!}v_{\ast}\overline{{\Bbb
Q}}_{\ell}[1]$ est canoniquement isomorphe au complexe $\Psi
=\beta_{\ast}\overline{{\Bbb Q}}_{\ell}$ d\'{e}fini en $1.3$.
\endthm

\rem D\'{e}monstration
\endrem
Le complexe $\beta^{\ast}f_{!}v_{\ast}\overline{{\Bbb
Q}}_{\ell}$ est l'image directe \`{a} supports propres par le
morphisme structural ${\Bbb G}_{{\rm m},S}\rightarrow S$ du faisceau
constant $\overline{{\Bbb Q}}_{\ell}$ sur ${\Bbb G}_{{\rm m},S}-a(S)$
prolong\'{e} par image directe totale \`{a} ${\Bbb G}_{{\rm m},S}$
tout entier. On a donc bien $\beta^{\ast}f_{!}v_{\ast}\overline{{\Bbb
Q}}_{\ell}\cong \overline{{\Bbb Q}}_{\ell}[-1]$ d'apr\`{e}s
l'assertion (i) du lemme 1.4.

En raisonnant comme dans la d\'{e}monstration du lemme 2.3, on voit
qu'il nous reste \`{a} v\'{e}rifier que $0^{\ast}\mathop{\rm pr}
\nolimits_{\ast}w_{!}\overline{{\Bbb Q}}_{\ell}=(0)$ o\`{u} $w$ est
l'inclusion dans ${\Bbb A}_{S}^{1}\times_{S}{\Bbb G}_{{\rm m},S}$ de
l'ouvert compl\'{e}mentaire du ferm\'{e} d'\'{e}quation $xt=a$.

Consid\'{e}rons la compactification partielle
$$\diagram{
{\Bbb A}_{S}^{1}\times_{S}{\Bbb G}_{{\rm m},S}&\kern -2mm
\smash{\mathop{\lhook\joinrel\mathrel{\hbox to 6mm{\rightarrowfill}}
}\limits^{\scriptstyle u}}\kern -2mm&{\Bbb A}_{S}^{1}\times_{S}({\Bbb
P}_{S}^{1}-0(S))\cr
\llap{$\scriptstyle \mathop{\rm pr}$}\left\downarrow
\vbox to 5mm{}\right.\rlap{}&&\cr
{\Bbb A}_{S}^{1}&&\cr
\arrow(28,16)\dir(-3,-2)\length{13}
\put (27,10){\scriptstyle\mathop{\widetilde{\rm pr}}}}
$$
le ferm\'{e} compl\'{e}mentaire de ${\Bbb A}_{S}^{1}\times_{S}{\Bbb
G}_{{\rm m},S}$ dans ${\Bbb A}_{S}^{1}\times_{S}({\Bbb
P}_{S}^{1}-0(S))$ \'{e}tant ${\Bbb A}_{S}^{1}\times_{S}\infty (S)$.
Si $y=t^{-1}$ est la coordonn\'{e}e affine \`{a} l'infini de ${\Bbb
P}_{S}^{1}$, le ferm\'{e} d'\'{e}quation $xt=a$ de ${\Bbb
A}_{S}^{1}\times_{S}{\Bbb G}_{{\rm m},S}$ se prolonge en le ferm\'{e}
d'\'{e}quation $x=ay$ de ${\Bbb A}_{S}^{1}\times_{S}({\Bbb
P}_{S}^{1}-0(S))=S[x,y]$, ferm\'{e} qui est transverse au ferm\'{e}
${\Bbb A}_{S}^{1}\times_{S}\infty (S)$ d'\'{e}quation $y=0$.  Si on
note $\widetilde{w}$ l'inclusion de l'ouvert $\{x\not=ay\}$ dans
${\Bbb A}_{S}^{1}\times_{S}({\Bbb P}_{S}^{1}-0(S))$ et $\widetilde{u}$
l'inclusion de $\{(x,y)\in {\Bbb A}_{S}^{1}\times_{S}{\Bbb G}_{{\rm
m},S}\mid x\not=ay\}$ dans $\{(x,y)\in {\Bbb
A}_{S}^{1}\times_{S}({\Bbb P}_{S}^{1}-0(S))\mid x\not=ay\}$, on a
$$
u_{\ast}w_{!}=\widetilde{w}_{!}\widetilde{u}_{\ast}\overline{{\Bbb
Q}}_{\ell}
$$
par la formule de K\"{u}nneth, et on a donc un triangle distingu\'{e}
$$
\mathop{\widetilde{\rm pr}}\nolimits_{\ast}\widetilde{w}_{!}
\overline{{\Bbb Q}}_{\ell}\rightarrow \mathop{\rm pr}
\nolimits_{\ast}w_{!}\overline{{\Bbb Q}}_{\ell}\rightarrow
\mathop{\widetilde{\rm pr}}\nolimits_{\ast}\widetilde{w}_{!}
R^{1}\widetilde{u}_{\ast}\overline{{\Bbb Q}}_{\ell}[-1]\rightarrow
$$

Or, on a
$$
\mathop{\widetilde{\rm pr}}\nolimits_{\ast}\widetilde{w}_{!}
\overline{{\Bbb Q}}_{\ell}=(0)
$$
d'apr\`{e}s 1.4(i), et $R^{1}\widetilde{u}_{\ast}\overline{{\Bbb
Q}}_{\ell}$ est la faisceau $\overline{{\Bbb Q}}_{\ell}(-1)$ sur la
partie localement ferm\'{e}e ${\Bbb G}_{{\rm m},S}\times_{S}\infty
(S)= \{(x,y)\in {\Bbb A}_{S}^{1}\times_{S}({\Bbb P}_{S}^{1}-0(S))\mid
x\not=0,~y=0\}$, de sorte que $\mathop{\widetilde{\rm pr}}
\nolimits_{\ast}\widetilde{w}_{!}
R^{1}\widetilde{u}_{\ast}\overline{{\Bbb Q}}_{\ell}[-1]$ est le
prolongement par z\'{e}ro de $\overline{{\Bbb Q}}_{\ell}[-1](-1)$ de
${\Bbb G}_{{\rm m},S}$ \`{a} ${\Bbb A}_{S}^{1}$.  Par suite, on a bien
$0^{\ast}\mathop{\rm pr}\nolimits_{\ast}w_{!}\overline{{\Bbb
Q}}_{\ell}=(0)$ et le lemme est d\'{e}montr\'{e}.
\hfill\hfill$\square$
\vskip 3mm

\thm LEMME 3.3
\enonce
On a
$$
\overline{h}_{!}\Psi =(0)
$$
o\`{u} $\overline{h}:{\cal A}_{S}\rightarrow S$ est le morphisme
structural.
\endthm

\rem D\'{e}monstration
\endrem
Cela r\'{e}sulte imm\'{e}diatement des lemmes 3.2 et 1.4(i) puisque
$h=\overline{h}\circ f$.

\hfill\hfill$\square$
\vskip 3mm

\thm LEMME 3.4
\enonce
Faisons agir ${\Bbb G}_{{\rm m},S}$ diagonalement par homoth\'{e}tie
sur ${\Bbb A}_{S}^{2}$ et soient $q:[{\Bbb A}_{S}^{2}/{\Bbb G}_{{\rm
m},S}]\rightarrow {\cal A}_{S}^{2}={\cal A}_{S}\times_{S}{\cal A}_{S}$
le morphisme quotient et $\sigma :[{\Bbb A}_{S}^{2}/{\Bbb G}_{{\rm
m},S}]\rightarrow {\cal A}_{S}$ le morphisme induit par le morphisme
diff\'{e}rence ${\Bbb A}_{S}^{2}\rightarrow {\Bbb A}_{S}^{1},~(x,y)
\mapsto x-y$.  Alors on a un isomorphisme canonique
$$
q_{!}\sigma^{\ast}\Psi[1]\cong\Psi\boxtimes\Psi.
$$
\endthm

\rem D\'{e}monstration
\endrem
On a d'une part
$$
\sigma^{\ast}\beta_{\ast}\overline{{\Bbb Q}}_{\ell}=j_{\ast}
\overline{{\Bbb Q}}_{\ell}
$$
o\`{u} $j:{\cal U}\hookrightarrow [{\Bbb A}_{S}^{2}/{\Bbb G}_{{\rm
m},S}]$ est l'ouvert quotient de l'ouvert $\{(x,y)\mid x\not=y\}$ de
${\Bbb A}_{S}^{2}$, puisque $\sigma$ est repr\'{e}sentable et lisse.
On a d'autre part
$$
\beta_{\ast}\overline{{\Bbb Q}}_{\ell}\boxtimes
\beta_{\ast}\overline{{\Bbb Q}}_{\ell}\cong
(\beta^{2})_{\ast}\overline{{\Bbb Q}}_{\ell}
$$
o\`{u} $\beta^{2}$ est l'immersion ouverte $\beta\times\beta: S=
S\times_{S}S\hookrightarrow {\cal A}_{S}^{2}$, d'apr\`{e}s la formule
de K\"{u}nneth.  On cherche donc \`{a} construire un isomorphisme
canonique
$$
q_{!}j_{\ast}\overline{{\Bbb Q}}_{\ell}[1]\buildrel\sim\over
\longrightarrow (\beta^{2})_{\ast}\overline{{\Bbb Q}}_{\ell}
$$

L'immersion ouverte $\beta^{2}$ se factorise en
$$
S\,\smash{\mathop{\lhook\joinrel\mathrel{\hbox to
8mm{\rightarrowfill}}}\limits^{\scriptstyle\gamma}}\,{\cal
B}\,\smash{\mathop{\lhook\joinrel\mathrel{\hbox to
8mm{\rightarrowfill}}}\limits^{\scriptstyle i}}\,{\cal A}_{S}^{2}
$$
o\`{u} $i:{\cal B}\hookrightarrow {\cal A}_{S}^{2}$ est l'image de
l'ouvert ${\Bbb A}_{S}^{2}-\{(0,0)\}\subset {\Bbb A}_{S}^{2}$.  Il 
suffit donc :
\vskip 1mm

\itemitem{(1)} de construire une isomorphisme canonique
$i^{\ast}q_{!}j_{\ast}\overline{{\Bbb Q}}_{\ell}[1]\buildrel\sim\over
\longrightarrow\gamma_{\ast}\overline{{\Bbb Q}}_{\ell}$,
\vskip 1mm

\itemitem{(2)} de montrer que la fl\`{e}che d'adjonction
$q_{!}j_{\ast}\overline{{\Bbb Q}}_{\ell}\rightarrow
i_{\ast}i^{\ast}q_{!}j_{\ast}\overline{{\Bbb Q}}_{\ell}$ est un
isomorphisme.
\vskip 1mm

Le champ ${\cal B}$ admet le recouvrement ouvert ${\cal B}= {\cal
B}_{1}\cup {\cal B}_{2}$ o\`{u} ${\cal B}_{1}={\cal A}_{S}\times_{S}S$
et ${\cal B}_{2}=S\times_{S}{\cal A}_{S}$ sont tous les deux
isomorphes \`{a} ${\cal A}_{S}$, et la restriction du morphisme compos\'{e}
$$
{\cal U}\,\smash{\mathop{\lhook\joinrel\mathrel{\hbox to
8mm{\rightarrowfill}}}\limits^{\scriptstyle j}}\, [{\Bbb
A}_{S}^{2}/{\Bbb G}_{{\rm m},S}]\,\smash{\mathop{\hbox to
8mm{\rightarrowfill}} \limits^{\scriptstyle q}}\,{\cal A}_{S}^{2}
$$
\`{a} l'un ou l'autre de ces ouverts est naturellement isomorphe \`{a}
$$
{\Bbb A}_{S}^{1}-\{1(S)\}\,\smash{\mathop{\lhook\joinrel\mathrel{\hbox
to 8mm{\rightarrowfill}} }\limits^{\scriptstyle v}}\,{\Bbb A}_{S}^{1}
\,\smash{\mathop{\hbox to 8mm{\rightarrowfill}}\limits^{\scriptstyle
f}}\,{\cal A}_{S}
$$
o\`{u} $1:S\rightarrow {\Bbb A}_{S}^{1}$ est la section constante de
valeur $1$.  Par suite, la restriction de $i^{\ast}q_{!}j_{\ast}
\overline{{\Bbb Q}}_{\ell}[1]$ \`{a} ${\cal B}_{1}\cong {\cal A}_{S}$
ou ${\cal B}_{2}\cong {\cal A}_{S}$ est canoniquement isomorphe \`{a}
$f_{\ast}v_{!}\overline{{\Bbb Q}}_{\ell}[1]$, et donc \`{a}
$\beta_{\ast}\overline{{\Bbb Q}}_{\ell}$, d'apr\`{e}s le lemme 3.2,
d'o\`{u} le point (1).

Il suffit de d\'{e}montrer le point (2) apr\`{e}s le changement de
base par le morphisme quotient ${\Bbb A}_{S}^{2}\rightarrow {\cal
A}_{S}^{2}$.  Par suite, d'apr\`{e}s le th\'{e}or\`{e}me de changement
de base par un morphisme lisse, si on note encore $j$ l'inclusion de
l'ouvert $U=\{(x,y,t)\mid xt\not=y\}$ dans ${\Bbb A}_{S}^{2}
\times_{S}{\Bbb G}_{{\rm m},S}$, $q$ la projection canonique de ${\Bbb
A}_{S}^{2}\times_{S}{\Bbb G}_{{\rm m},S}$ sur ${\Bbb A}_{S}^{2}$ et
$i$ l'inclusion de ${\Bbb A}_{S}^{2}-\{(0,0)\}$ dans ${\Bbb
A}_{S}^{2}$, il ne nous reste plus qu'\`{a} d\'{e}montrer que la
fl\`{e}che d'adjonction $q_{!}j_{\ast}\overline{{\Bbb
Q}}_{\ell}\rightarrow i_{\ast}i^{\ast}q_{!}j_{\ast}\overline{{\Bbb
Q}}_{\ell}$ est un isomorphisme dans $D_{{\rm c}}^{{\rm b}}({\Bbb
A}_{S}^{2},\overline{{\Bbb Q}}_{\ell})$.

On compactifie la projection $q$ en
$$\diagram{
{\Bbb A}_{S}^{2}\times_{S}{\Bbb G}_{{\rm
m},S}&\kern -2mm
\smash{\mathop{\lhook\joinrel\mathrel{\hbox to 6mm{\rightarrowfill}}
}\limits^{\scriptstyle k}}\kern -2mm&{\Bbb A}_{S}^{2}\times_{S}{\Bbb
P}_{S}^{1}\cr
\llap{$\scriptstyle q$}\left\downarrow
\vbox to 5mm{}\right.\rlap{}&&\cr
{\Bbb A}_{S}^{2}&&\cr
\arrow(27,16)\dir(-3,-2)\length{13}
\put (26,10.7){\scriptstyle \overline{q}}}
$$
et on introduit l'ouvert
$$
\overline{j}:\overline{U}=\{((x,y),(T;U))\mid xT-yV\not=0\}
\hookrightarrow {\Bbb A}_{S}^{2}\times_{S}{\Bbb P}_{S}^{1}.
$$

On a un carr\'{e} cart\'{e}sien
d'immersions ouvertes
$$\diagram{
U&\kern -1mm\smash{\mathop{\hbox to 8mm{\rightarrowfill}}
\limits^{\scriptstyle j}}\kern -1mm&{\Bbb A}_{S}^{2}\times_{S}
{\Bbb G}_{{\rm m},S}\cr
\llap{$\scriptstyle k_{U}$}\left\downarrow
\vbox to 4mm{}\right.\rlap{}&\square&\llap{}\left\downarrow
\vbox to 4mm{}\right.\rlap{$\scriptstyle k$}\cr
\overline{U}&\kern -1mm\smash{\mathop{\hbox to 8mm{\rightarrowfill}}
\limits_{\scriptstyle\overline{j}}}\kern -1mm&{\Bbb A}_{S}^{2}
\times_{S}{\Bbb P}_{S}^{1}\cr}
$$
et, par adjonction, une fl\`{e}che
$$
k_{!}j_{\ast}\overline{{\Bbb Q}}_{\ell}\rightarrow\overline{j}_{\ast}
k_{U,!}\overline{{\Bbb Q}}_{\ell}.
$$

Cette fl\`{e}che est un isomorphisme.  En effet, le probl\`{e}me est
local sur ${\Bbb A}_{S}^{2} \times_{S}{\Bbb P}_{S}^{1}$ en les points
du ferm\'{e} compl\'{e}mentaire de la r\'{e}union des deux ouverts
${\Bbb A}_{S}^{2} \times_{S}{\Bbb A}_{S}^{1}$ et $\overline{U}$.  Or,
en ces points on est localement dans une situation produit.  En effet,
on a des coordonn\'{e}es $(x,y,(T;U))$ sur ${\Bbb A}_{S}^{2}
\times_{S}{\Bbb P}_{S}^{1}$, ${\Bbb A}_{S}^{2} \times_{S}{\Bbb
A}_{S}^{1}$ est le compl\'{e}mentaire du diviseur lisse
$\{(T;U)=(0;1)\} \cup\{(T;U)=(1;0)\}$, $\overline{U}$ est le
compl\'{e}mentaire du diviseur lisse $\{xT-yU=0\}$ et ces deux
diviseurs lisses se coupent transversalement.  On conclut donc par la
formule de K\"{u}nneth.

On a donc des isomorphismes canoniques
$$
q_{!}j_{\ast}\overline{{\Bbb Q}}_{\ell}\cong\overline{q}_{\ast}k_{!}
j_{\ast}\overline{{\Bbb Q}}_{\ell}\cong\overline{q}_{\ast}
\overline{j}_{\ast} k_{U,!}\overline{{\Bbb Q}}_{\ell}
$$
et, comme $\overline{q}\circ\overline{j}$ se factorise par $i$ en
$$\diagram{
\overline{U}&\kern -10mm\smash{\mathop{\hbox to 17mm
{\rightarrowfill}}\limits^{\scriptstyle\overline{j}}}\kern
-1mm&{\Bbb A}_{S}^{2} \times_{S}{\Bbb P}_{S}^{1}\cr
\llap{$\scriptstyle \overline{q}_{\overline{U}}$}\left\downarrow
\vbox to 4mm{}\right.\rlap{}&&\llap{}\left\downarrow
\vbox to 4mm{}\right.\rlap{$\scriptstyle\overline{q}$}\cr
{\Bbb A}_{S}^{2}-\{(0,0)\}&\kern -1mm\smash{\mathop{\hbox to
14mm{\rightarrowfill}} \limits_{\scriptstyle i}}\kern -7mm&{\Bbb
A}_{S}^{2}\cr}
$$
on a en fait un isomorphisme canonique
$$
q_{!}j_{\ast}\overline{{\Bbb Q}}_{\ell}\cong i_{\ast}
\overline{q}_{\overline{U},\ast} k_{U,!}\overline{{\Bbb Q}}_{\ell},
$$
d'o\`{u} le point (2).
\hfill\hfill$\square$
\vskip 3mm

\rem D\'{e}monstration du th\'{e}or\`{e}me $3.1$
\endrem
On proc\`{e}de comme pour la transformation de Fourier-Deligne
(cf. [8](1.2.2)).  Il s'agit donc de construire un isomorphisme
$$
\mathop{\rm pr}\nolimits_{13,!}\nu^{\ast}(\Psi\boxtimes\Psi )\cong
\Delta_{!}\overline{{\Bbb Q}}_{\ell}[2-2r](-r)
$$
o\`{u} $\mathop{\rm pr}\nolimits_{13}:{\cal V}\times_{S}{\cal
V}^{\vee}\times_{S}{\cal V}\rightarrow {\cal V}\times_{S}{\cal V}$ est
le projection canonique, o\`{u}
$$
\nu :{\cal V}\times_{S}{\cal V}^{\vee}\times_{S}{\cal V}\rightarrow
{\cal A}_{S}^{2}
$$
est induit par le morphisme $(v_{1},w,v_{2})\rightarrow (\langle
w,v_{1}\rangle ,\langle w,v_{2}\rangle )$ de $V\times_{S} V^{\vee}
\times_{S}V$ dans ${\Bbb A}_{S}^{2}$ et o\`{u} $\Delta :{\cal V}
\rightarrow {\cal V}\times_{S}{\cal V}$ est le morphisme diagonal.

On a le carr\'{e} cart\'{e}sien
$$\diagram{
{\cal V}^{\vee}\times_{S}[(V\times_{S}V)/{\Bbb G}_{{\rm m},S}]&\kern
-1mm\smash{\mathop{\hbox to 8mm{\rightarrowfill}}
\limits^{\scriptstyle\nu'}}\kern -1mm& [{\Bbb A}_{S}^{2}/{\Bbb
G}_{{\rm m},S}]\cr
\llap{$\scriptstyle q'$}\left\downarrow
\vbox to 4mm{}\right.\rlap{}&\square&\llap{}\left\downarrow
\vbox to 4mm{}\right.\rlap{$\scriptstyle q$}\cr
{\cal V}\times_{S}{\cal V}^{\vee}\times_{S}{\cal V}&\kern -8mm
\smash{\mathop{\hbox to 20mm{\rightarrowfill}}\limits_{\scriptstyle
\nu}}\kern -6mm&{\cal A}_{S}^{2}\cr}
$$
o\`{u} $\nu'$ est induite par le morphisme de
$V^{\vee}\times_{S}V\times_{S}V$ dans ${\Bbb A}_{S}^{2}$ qui envoie
$(w,v_{1},v_{2})$ sur $(\langle w,v_{1}\rangle ,\langle w,v_{2}\rangle
)$ et $q'$ est induite par le morphisme de $V^{\vee}\times_{S}V
\times_{S} V$ dans $V\times_{S}V^{\vee}\times_{S}V$ qui envoie
$(w,v_{1},v_{2})$ sur $(v_{1},w,v_{2})$.  On d\'{e}duit de
l'isomorphisme du lemme 3.4 et du th\'{e}or\`{e}me de changement de
base propre un isomorphisme
$$
\nu^{\ast}(\Psi\boxtimes\Psi)\cong q_{!}'\nu'^{\ast}\sigma^{\ast}\Psi
[1].
$$
Comme $\sigma\circ\nu'$ se factorise en
$$
{\cal V}^{\vee}\times_{S}[(V\times_{S}V)/{\Bbb G}_{{\rm 
m},S}]\,\smash{\mathop{\hbox
to 12mm{\rightarrowfill}}\limits^{\scriptstyle\mathop{\rm Id}\times
s}}\, {\cal V}\times_{S}{\cal V}^{\vee}\,\smash{\mathop{\hbox to
8mm{\rightarrowfill}}\limits^{\scriptstyle\mu}}\,{\cal A}_{S}
$$
o\`{u} $s:[(V\times_{S}V)/{\Bbb G}_{{\rm m},S}]\rightarrow {\cal V}$ 
est induit par
le morphisme $V\times_{S}V\rightarrow V$ qui envoie $(v_{1},v_{2})$
sur $v_{1}-v_{2}$, que $\mathop{\rm pr}\nolimits_{13}\circ q'$ se
factorise en
$$
{\cal V}^{\vee}\times_{S}[(V\times_{S}V)/{\Bbb G}_{{\rm 
m},S}]\,\smash{\mathop{\hbox
to 8mm{\rightarrowfill}}\limits^{\scriptstyle\mathop{\rm pr}
\nolimits_{2}}}\, [(V\times_{S}V)/{\Bbb G}_{{\rm m},S}]\rightarrow {\cal V}
\times_{S}{\cal V}
$$
et que $\Delta$ se factorise en
$$
{\cal V}\,\smash{\mathop{\hbox to 8mm{\rightarrowfill}}
\limits^{\scriptstyle\Delta'}}\,[(V\times_{S}V)/{\Bbb G}_{{\rm
m},S}]\rightarrow {\cal V}\times_{S}{\cal V},
$$
o\`{u} $\Delta'$ est induit par le morphisme diagonal $V\rightarrow
V\times_{S}V$, il suffit de construire un isomorphisme
$$
\mathop{\rm pr}\nolimits_{2,!}(\mathop{\rm Id}\times s)^{\ast}\mu^{\ast}
\Psi [1]\cong\Delta_{!}'\overline{{\Bbb Q}}_{\ell}[2-2r](-r).
$$
Or, on a un diagramme \`{a} carr\'{e}s cart\'{e}siens
$$\diagram{
[(V\times_{S}V)/{\Bbb G}_{{\rm m},S}]\times_{S}{\cal V}^{\vee}&\kern -3mm
\smash{\mathop{\hbox to 14mm{\rightarrowfill}}\limits^{\scriptstyle
\mathop{\rm Id}\times s}}\kern -1mm&{\cal V}^{\vee}\times_{S}{\cal
V}\cr
\llap{$\scriptstyle\mathop{\rm pr}\nolimits_{2}$}\left\downarrow
\vbox to 4mm{}\right.\rlap{}&\square&\llap{}\left\downarrow\vbox to
4mm{}\right.\rlap{$\scriptstyle\mathop{\rm pr}\nolimits_{2}$}\cr
[(V\times_{S}V)/{\Bbb G}_{{\rm m},S}]&\kern -6mm\smash{\mathop{\hbox to
22mm{\rightarrowfill}}\limits^{\scriptstyle s}}\kern -5mm&{\cal V}\cr
\llap{$\scriptstyle\Delta'$}\left\uparrow\vbox to 4mm{}
\right.\rlap{}&\square&\llap{}\left\uparrow\vbox to
4mm{}\right.\rlap{$\scriptstyle i$}\cr
{\cal V}&\kern -17mm\smash{\mathop{\hbox to 30mm{\rightarrowfill}}
\limits_{\scriptstyle\pi}}\kern -3mm&B({\Bbb G}_{{\rm m},S})\cr}
$$
o\`{u} $i$ est induit par la section nulle $S\rightarrow V$. Par
suite, en utilisant de nouveau le th\'{e}or\`{e}me de changement de
base propre, on voit qu'il suffit de construire un isomorphisme
$$
\mathop{\rm pr}\nolimits_{2,!}\mu^{\ast}\Psi [1]\cong
i_{!}\overline{{\Bbb Q}}_{\ell}[2-2r](-r).
$$

Comme $i$ est repr\'{e}sentable et une immersion ferm\'{e}e, on
obtient un tel isomorphisme si l'on montre que $\mathop{\rm
pr}\nolimits_{2,!}\mu^{\ast}\Psi [1]$ est support\'{e} par l'image de
$i$ et que, apr\`{e}s changement de base par $i$, on a un tel
isomorphisme.

Or, pour tout point g\'{e}om\'{e}trique $s$ de $S$ et tout $v\in {\cal
V}_{s}$, $v\not=0$, induisant un morphisme $\langle v\rangle :{\cal
V}_{s}^{\vee}\rightarrow {\cal A}_{s},~w\rightarrow\langle
w,v\rangle$, on a
$$
R\Gamma_{{\rm c}}({\cal V}_{s}^{\vee},\langle v\rangle^{\ast}\Psi )
\cong R\Gamma_{{\rm c}}({\cal A}_{s},\langle v\rangle_{!}\langle
v\rangle^{\ast}\Psi )\cong R\Gamma_{{\rm c}}({\cal A}_{s},\Psi
)[2-2r](1-r)=(0)
$$
d'apr\`{e}s le lemme 3.3, et on a
$$
(\overline{\pi}^{\vee})_{!}(\pi^{\vee})^{\ast}\alpha^{\ast}\Psi [1]\cong
\varepsilon_{!}(\pi^{\vee})_{!}(\pi^{\vee})^{\ast}\alpha^{\ast}\Psi [1]\cong
\varepsilon_{!}\alpha^{\ast}\Psi [1-2r](-r)
$$
o\`{u} le morphisme structural $\overline{\pi}^{\vee}:{\cal
V}^{\vee}\rightarrow S$ est le compos\'{e} de $\pi^{\vee}:{\cal
V}^{\vee}\rightarrow B({\Bbb G}_{{\rm m},S})$ par le morphisme structural
$\varepsilon :B({\Bbb G}_{{\rm m},S})\rightarrow S$.  Mais d'apr\`{e}s le lemme
1.4(ii), on a un isomorphisme canonique
$$
\alpha^{\ast}\Psi\cong g_{!}\overline{{\Bbb Q}}_{\ell}[1]
$$
o\`{u} $g:S\rightarrow B({\Bbb G}_{{\rm m},S})$ est le torseur 
universel, d'o\`{u} la
conclusion puisque $\varepsilon\circ g$ est l'identit\'{e} de $S$.
\hfill\hfill$\square$
\vskip 3mm

\section{4}{Commutation \`{a} la dualit\'{e} et $t$-exactitude de la
transformation de Fourier homog\`{e}ne}

\thm TH\'{E}OR\`{E}ME 4.1
\enonce
Il existe un isomorphisme canonique entre les foncteurs
$$
K\mapsto\mathop{\rm Four}\nolimits_{{\cal V}/S}(D_{{\cal V}}(K))
\hbox{ et }K\mapsto D_{{\cal V}^{\vee}}(\mathop{\rm Four}
\nolimits_{{\cal V}/S}(K))(-r)
$$
de $D_{{\rm c}}^{{\rm b}}({\cal V},\overline{{\Bbb Q}}_{\ell})^{{\rm
opp}}$ dans $D_{{\rm c}}^{{\rm b}}({\cal V}^{\vee},\overline{{\Bbb
Q}}_{\ell})$.
\endthm

\rem D\'{e}monstration
\endrem
On reprend mot pour mot un argument de Verdier pour la transformation
de Fourier-Deligne. Par la dualit\'{e} de Grothendieck, le foncteur
$$
D_{{\cal V}^{\vee}}\circ\mathop{\rm Four}\nolimits_{{\cal V}/S}\circ
D_{{\cal V}}:D_{{\rm c}}^{{\rm b}}({\cal V},\overline{{\Bbb Q}}_{\ell})
\rightarrow D_{{\rm c}}^{{\rm b}}({\cal V}^{\vee},\overline{{\Bbb
Q}}_{\ell})
$$
est isomorphe au foncteur ${\cal F}_{\ast}$ d\'{e}fini par
$$
{\cal F}_{\ast}(K)=\mathop{\rm pr}\nolimits_{\ast}^{\vee} (\mathop{\rm
pr}\nolimits^{!}K\widetilde{\otimes}\mu^{!}\Psi')[r-1]
$$
o\`{u} $\Psi'=D_{{\cal A}}(\Psi )=\beta_{!}\overline{{\Bbb Q}}_{\ell}$.

Or on a un isomorphisme canonique bi-fonctoriel en $K$ dans $D_{{\rm
c}}^{{\rm b}}({\cal V},\overline{{\Bbb Q}}_{\ell})$ et $L$ dans
$D_{{\rm c}}^{{\rm b}}({\cal V}^{\vee},\overline{{\Bbb
Q}}_{\ell})^{{\rm opp}}$,
$$
\mathop{\rm Hom}\nolimits_{D_{{\rm c}}^{{\rm b}}({\cal V},
\overline{{\Bbb Q}}_{\ell})}(\mathop{\rm Four}\nolimits_{{\cal
V}^{\vee}/S}(L),K)=\mathop{\rm Hom}\nolimits_{D_{{\rm c}}^{{\rm
b}}({\cal V}^{\vee},\overline{{\Bbb Q}}_{\ell})}(L,{\cal F}_{\ast}(K))
$$
puisque $\mathop{\rm pr}\nolimits_{!}$ est l'adjoint \`{a} gauche de
$\mathop{\rm pr}\nolimits^{!}$ et que l'on a l'isomorphisme de
Grothendieck {\og}{cher \`{a} Cartan}{\fg}
$$\eqalign{
\mathop{\rm Hom}\nolimits_{D_{{\rm c}}^{{\rm b}}({\cal X},
\overline{{\Bbb Q}}_{\ell})}(K\otimes L,M)&\cong\mathop{\rm
Hom}\nolimits_{D_{{\rm c}}^{{\rm b}}({\cal X},\overline{{\Bbb
Q}}_{\ell})}(K,\mathop{{\cal H}{\it om}}(L,M))\cr
&\cong\mathop{\rm Hom}\nolimits_{D_{{\rm c}}^{{\rm b}}({\cal X},
\overline{{\Bbb Q}}_{\ell})}(K,D_{{\cal X}}(L)
\widetilde{\otimes}M)\cr}
$$
pour tout $S$-champ alg\'{e}brique de type fini et tous $K,L,M\in
\mathop{\rm ob}D_{{\rm c}}^{{\rm b}}({\cal X},\overline{{\Bbb
Q}}_{\ell})$.

Par suite, ${\cal F}_{\ast}$ est un adjoint \`{a} droite du foncteur
$\mathop{\rm Four}\nolimits_{{\cal V}^{\vee}/S}$.  Mais, il en est de
m\^{e}me du foncteur $K\rightarrow\mathop{\rm Four}\nolimits_{{\cal
V}/S}(K)(r)$ qui est un quasi-inverse de $\mathop{\rm
Four}\nolimits_{{\cal V}^{\vee}/S}$.  Ces deux adjoints \`{a} droite
sont donc quasi-isomorphes et le th\'{e}or\`{e}me est d\'{e}montr\'{e}.
\hfill\hfill$\square$
\vskip 3mm

\thm TH\'{E}OR\`{E}ME 4.2
\enonce
La transformation de Fourier homog\`{e}ne $\mathop{\rm
Four}\nolimits_{{\cal V}/S}$ est $t$-exacte.
\endthm

\rem D\'{e}monstration
\endrem
Compte tenu du th\'{e}or\`{e}me 4.1, il suffit de d\'{e}montrer que,
pour tout $\overline{{\Bbb Q}}_{\ell}$-faisceau pervers $K$ sur ${\cal
V}$, $\mathop{\rm Four}\nolimits_{{\cal V}/S}(K)$ est dans ${}^{{\rm
p}}D_{{\rm c}}^{\geq 0}({\cal V}^{\vee},\overline{{\Bbb Q}}_{\ell})$,
c'est-\`{a}-dire que les faisceaux de cohomologie perverse ${}^{{\rm
p}}{\cal H}^{i}(\mathop{\rm Four}\nolimits_{{\cal V}/S}(K))$ sont tous
nuls pour $i<0$.

Posons $L=\mathop{\rm pr}\nolimits^{\ast}K\otimes\mu^{\ast}\Psi$ et
d\'{e}montrons dans un premier temps que $L$ est dans ${}^{{\rm
p}}D_{{\rm c}}^{\geq 0}({\cal V}^{\vee}\times_{S}{\cal
V},\overline{{\Bbb Q}}_{\ell})$.  On a le d\'{e}vissage
$$
\overline{{\Bbb Q}}_{\ell}\rightarrow\beta_{\ast}\overline{{\Bbb
Q}}_{\ell}\rightarrow\alpha_{\ast}\overline{{\Bbb Q}}_{\ell}[-1](-1)
\rightarrow
$$
o\`{u} on rappelle que $\alpha :B({\Bbb G}_{{\rm m},S})\hookrightarrow
{\cal A}_{S}$ est l'immersion ferm\'{e}e compl\'{e}mentaire de l'immersion
ouverte $\beta :S\hookrightarrow {\cal A}$, et on a donc aussi le
d\'{e}vissage
$$
\mathop{\rm pr}\nolimits^{\ast}K[r-1]\rightarrow L\rightarrow
i_{\ast}i^{\ast}(\mathop{\rm pr}\nolimits^{\ast}K[r-1])[-1](-1)
\rightarrow
$$
o\`{u} $i:{\cal H}\hookrightarrow {\cal V}^{\vee}\times_{S}{\cal V}$
est l'immersion ferm\'{e}e induite par l'inclusion $\{\langle
w,v\rangle =0\}\subset V^{\vee}\times_{S}V$.  Maintenant, comme $K$
est pervers et que $\mathop{\rm pr}$ est lisse purement de dimension
relative $r-1$, $\mathop{\rm pr}\nolimits^{\ast}K[r-1]$ est pervers et
$i_{\ast}i^{\ast}(\mathop{\rm pr}\nolimits^{\ast}K[r-1])$ est dans
${}^{{\rm p}}D_{{\rm c}}^{[-1,0]}({\cal V}^{\vee}\times_{S}{\cal
V},\overline{{\Bbb Q}}_{\ell})$ d'apr\`{e}s [6](4.1.10)(ii).  Par
suite, $L$ est bien dans ${}^{{\rm p}}D_{{\rm c}}^{\geq 0}({\cal
V}^{\vee}\times_{S}{\cal V},\overline{{\Bbb Q}}_{\ell})$.

Pour terminer la d\'{e}monstration, \'{e}tudions le complexe
$\pi_{!}\rho^{\ast}L[1]$ o\`{u} on a not\'{e} encore $\pi :{\cal
V}^{\vee}\times_{S}V\rightarrow {\cal V}^{\vee}$ et $\rho : {\cal
V}^{\vee}\times_{S} V\rightarrow {\cal V}^{\vee}\times_{S}{\cal V}$
les changements de base par $\overline{\pi}^{\vee}:{\cal
V}^{\vee}\rightarrow S$ des morphismes $\pi :V\rightarrow S$ et $\rho
:V\rightarrow {\cal V}$.  Comme $\rho$ est lisse purement de dimension
relative $1$, $\rho^{\ast}L[1]$ est pervers sur ${\cal
V}^{\vee}\times_{S}V$ et, comme $\pi$ est affine,
$\pi_{!}\rho^{\ast}L[1]$ est donc dans ${}^{{\rm p}} D_{{\rm c}}^{\geq
0}({\cal V}^{\vee},\overline{{\Bbb Q}}_{\ell})$ d'apr\`{e}s [6](4.1.2).
De plus, on a $\pi_{!}\rho^{\ast}L[1]=\mathop{\rm pr}
\nolimits_{!}(L\otimes\rho_{!}\overline{{\Bbb Q}}_{\ell}[1])$
d'apr\`{e}s la formule des projections et on a un triangle
distingu\'{e}
$$
\overline{{\Bbb Q}}_{\ell}\rightarrow\rho_{!}\overline{{\Bbb
Q}}_{\ell} [1]\rightarrow\overline{{\Bbb Q}}_{\ell}[-1](-1)\rightarrow .
$$
Par cons\'{e}quent, on a un triangle distingu\'{e}
$$
\mathop{\rm pr}\nolimits_{!}^{\vee}L\rightarrow\pi_{!}\rho^{\ast}L[1]
\rightarrow\mathop{\rm pr}\nolimits_{!}^{\vee}L[-1](-1)\rightarrow
$$
o\`{u} $\pi_{!}\rho^{\ast}L[1]$ est dans ${}^{{\rm p}} D_{{\rm
c}}^{\geq 0}({\cal V}^{\vee},\overline{{\Bbb Q}}_{\ell})$.  En prenant
la suite exacte longue de cohomologie perverse, on en d\'{e}duit des
isomorphismes
$$
{}^{{\rm p}}{\cal H}^{i-2}\mathop{\rm pr}\nolimits_{!}^{\vee}L(-1)
\buildrel\sim\over\longrightarrow {}^{{\rm p}}{\cal H}^{i}\mathop{\rm
pr}\nolimits_{!}^{\vee}L,~\forall i<0.
$$
Par suite, si l'un des faisceaux de cohomologie perverse ${}^{{\rm
p}}{\cal H}^{i}$ de $\mathop{\rm pr}\nolimits_{!}^{\vee}L=\mathop{\rm
Four}\nolimits_{{\cal V}/S}(K)$ avec $i<0$ \'{e}tait non nul, il en
serait de m\^{e}me de tous les ${}^{{\rm p}}{\cal H}^{i-2n}$ avec
$n\geq 0$, ce qui contredit le fait que $\mathop{\rm Four}
\nolimits_{{\cal V}/S}(K)$ est un complexe born\'{e} (cf.  Corollaire
1.9), d'o\`{u} la conclusion.
\hfill\hfill$\square$

\section{5}{Compl\'{e}ments}

5.1. Supposons tout d'abord que $S=\mathop{\rm Spec}({\Bbb F}_{q})$
pour un corps fini \`{a} $q$ \'{e}l\'{e}ments ${\Bbb F}_{q}$ et
calculons l'action de la transformation de Fourier homog\`{e}ne sur
les fonctions {\og}{trace de Frobenius}{\fg}.

On note $|{\cal V}({\Bbb F}_{q})|$ l'ensemble des classes d'isomorphie
des objets de la cat\'{e}gorie ${\cal V}({\Bbb F}_{q})$; on a
$$
|{\cal V}({\Bbb F}_{q})|=\{0\}\cup {\Bbb P}(V)({\Bbb F}_{q}).
$$
Pour tout $K\in D_{{\rm c}}^{{\rm b}}({\cal V},\overline{{\Bbb
Q}}_{\ell})$, on note
$$
\tau_{K}:|{\cal V}({\Bbb F}_{q})|\rightarrow\overline{{\Bbb
Q}}_{\ell}
$$
la fonction qui associe \`{a} $v\in |{\cal V}({\Bbb F}_{q})|$ la trace
de l'action de l'endomorphisme de Frobenius g\'{e}om\'{e}trique
$\mathop{\rm Frob}\nolimits_{q}$ relatif \`{a} ${\Bbb F}_{q}$ sur la
fibre de $K$ en n'importe quel point g\'{e}om\'{e}trique de ${\cal V}$
localis\'{e} en $v$.  Bien entendu on pose la m\^{e}me d\'{e}finition
pour ${\cal V}^{\vee}$.

\thm LEMME 5.2
\enonce
Pour tout $K\in D_{{\rm c}}^{{\rm b}}({\cal V},\overline{{\Bbb
Q}}_{\ell})$ et tout $w\in |{\cal V}^{\vee}({\Bbb F}_{q})|$, on a
$$
\tau_{\mathop{\rm Four}\nolimits_{{\cal V}/S}(K)}(w)=\tau_{K}(0)-
\sum_{v\in {\Bbb P}(V)({\Bbb F}_{q})}\tau_{K}(v)+q
\sum_{{\scriptstyle v\in {\Bbb P}(V)({\Bbb F}_{q}) \atop\scriptstyle
(w,v)\in H({\Bbb F}_{q})}}\tau_{K}(v)
$$
o\`{u} on rappelle que $H\subset {\Bbb P}(V^{\vee})\times_{S}{\Bbb
P}(V)$ est l'hypersurface d'incidence.
\endthm

\rem Preuve
\endrem
On applique la formule des traces de Grothendieck, soit directement,
soit \`{a} travers la compatibilit\'{e} de la transformation de
Fourier homog\`{e}ne \`{a} celle de Fourier-Deligne prouv\'{e}e dans
la section 2.
\hfill\hfill$\square$
\vskip 3mm

5.3.  On revient \`{a} la situation g\'{e}n\'{e}rale et on
consid\`{e}re la sous-cat\'{e}gorie pleine
$$
D_{{\rm c},\dagger}^{{\rm b}}({\Bbb P}(V),\overline{{\Bbb Q}}_{\ell})
\subset D_{{\rm c}}^{{\rm b}}({\Bbb P}(V),\overline{{\Bbb Q}}_{\ell})
$$
d\'{e}finie par la condition $(\overline{\pi}^{\circ})_{!}K=
(\overline{\pi}^{\circ})_{\ast}K=(0)$, o\`{u} on rappelle que le
morphisme projectif $\overline{\pi}^{\circ}: {\cal V}^{\circ}={\Bbb
P}(V)\rightarrow S$ est la projection canonique.  Cette
sous-cat\'{e}gorie est triangul\'{e}e et est sa propre image par la
dualit\'{e} $D_{{\Bbb P}(V)}$.

La proposition suivante est une variante de r\'{e}sultats de
Brylinski sur la transformation de Radon g\'{e}om\'{e}trique (cf.
[2]).

\thm PROPOSITION 5.4
\enonce
{\rm (i)} Pour tout $K\in\mathop{\rm ob}D_{{\rm c}}^{{\rm b}} ({\Bbb
P}(V),\overline{{\Bbb Q}}_{\ell})$, on a
$(\overline{\pi}^{\circ})_{\ast}K=(0)$ si et seulement si la
fl\`{e}che d'oubli des supports
$$
j_{!}K\rightarrow j_{\ast}K
$$
est un isomorphisme dans $D_{{\rm c}}^{{\rm b}}({\cal V},
\overline{{\Bbb Q}}_{\ell})$.

\decale{\rm (ii)} Pour tout $K\in\mathop{\rm ob}D_{{\rm c}}^{{\rm b}}
({\Bbb P}(V),\overline{{\Bbb Q}}_{\ell})$, on a
$(\overline{\pi}^{\circ})_{\ast}K=(0)$ si et seulement si
$(\overline{\pi}^{\vee\circ})_{\ast}\mathop{\rm Rad}\nolimits_{{\Bbb
P}(V)/S}(K)=(0)$, et la transformation de Radon g\'{e}om\'{e}trique
$\mathop{\rm Rad}\nolimits_{{\Bbb P}(V)/S}$ induit une \'{e}quivalence
de cat\'{e}gories triangul\'{e}es
$$
\mathop{\rm Rad}\nolimits_{{\Bbb P}(V^{\vee})/S,\dagger}:D_{{\rm c},
\dagger}^{{\rm b}}({\Bbb P}(V),\overline{{\Bbb Q}}_{\ell})\buildrel
\sim\over\longrightarrow D_{{\rm c},\dagger}^{{\rm b}}({\Bbb
P}(V^{\vee}),\overline{{\Bbb Q}}_{\ell}),
$$
de quasi-inverse induit par $\mathop{\rm Rad}\nolimits_{{\Bbb
P}(V^{\vee})/S}(r-2)$.
\endthm

Nous utiliserons le lemme d'homotopie bien connu suivant:

\thm LEMME 5.5
\enonce
Soit $f:X\rightarrow S$ un $S$-sch\'{e}ma s\'{e}par\'{e} et de type
fini, muni d'une $S$-action contractante de ${\Bbb G}_{{\rm m},S}$,
c'est-\`{a}-dire d'un morphisme ${\Bbb A}_{S}^{1}\times_{S}
X\rightarrow X,~(u,x)\mapsto u\cdot x$ tel que $1\cdot x=x$ et
$(tu)\cdot x=t\cdot (u\cdot x)$ quels que soient $x\in X$, $t\in {\Bbb
G}_{{\rm m},S}$ et $u\in {\Bbb A}_{S}^{1}$.  Notons $\rho :
X\rightarrow {\cal X}$ le champ quotient de $X$ par cette action de
${\Bbb G}_{{\rm m},S}$.  Alors, pour tout $K\in \mathop{\rm ob}D_{{\rm
c}}^{{\rm b}}({\cal X},\overline{{\Bbb Q}}_{\ell})$ la
fl\`{e}che de restriction
$$
f_{\ast}\rho^{\ast}K\rightarrow g_{\ast}(\rho^{\ast}K|Y)
$$
au ferm\'{e}
$$
g:Y=\{x\in X\mid t\cdot u=x,~\forall t\in {\Bbb G}_{{\rm
m},S}\}\rightarrow S
$$
des points fixes de cette action, est un isomorphisme.
\endthm

\rem D\'{e}monstration
\endrem
Par d\'{e}vissage de $K$ en un complexe nul sur ${\cal Y}=[Y/{\Bbb
G}_{{\rm m},S}]\subset {\cal X}$ et un complexe support\'{e} sur ${\cal
Y}$, on peut supposer que la restriction de $K$ \`{a} ${\cal Y}$ est
nulle, et il s'agit alors de d\'{e}montrer que
$f_{\ast}\rho^{\ast}K=(0)$.

Soit $\tau :{\Bbb A}_{S}^{1}\times_{S}X\rightarrow {\Bbb
A}_{S}^{1}\times_{S}X$ l'endomorphisme d\'{e}fini par $\tau
(u,x)=(u,u\cdot x)$.  On a un isomorphisme \'{e}vident entre les
restrictions des complexes $\mathop{\rm pr}\nolimits_{X}^{\ast}
\rho^{\ast}K$ et $\tau^{\ast} \mathop{\rm pr}\nolimits_{X}^{\ast}
\rho^{\ast}K$ \`{a} l'ouvert ${\Bbb G}_{{\rm m},S}\times_{S}X\subset
{\Bbb A}_{S}^{1} \times_{S}X$, et comme la restriction de $K$ \`{a}
${\cal Y}$ est nulle par hypoth\`{e}se, il en est de m\^{e}me de la
restriction de $\tau^{\ast} \mathop{\rm pr}\nolimits_{X}^{\ast}
\rho^{\ast}K$ au ferm\'{e} $X\hookrightarrow {\Bbb
A}_{S}^{1}\times_{S}X,~x\mapsto (0,x)$.  On en d\'{e}duit l'existence
d'une fl\`{e}che
$$
\tau^{\ast} \mathop{\rm pr}\nolimits_{X}^{\ast}\rho^{\ast}K
\rightarrow\mathop{\rm pr}\nolimits_{X}^{\ast}\rho^{\ast}K
$$
sur ${\Bbb A}_{S}^{1}\times_{S}X$ tout entier dont la restriction
\`{a} ${\Bbb G}_{{\rm m},S}\times_{S}X$ est un isomorphisme. Par
adjonction, cette fl\`{e}che nous donne une fl\`{e}che
$$
\mathop{\rm pr}\nolimits_{X}^{\ast}\rho^{\ast}K\rightarrow \tau_{\ast}
\mathop{\rm pr}\nolimits_{X}^{\ast}\rho^{\ast}K
$$
et donc une fl\`{e}che
$$
(\mathop{\rm pr}\nolimits_{{\Bbb A}_{S}^{1}})_{\ast}
\mathop{\rm pr}\nolimits_{X}^{\ast}\rho^{\ast}K\rightarrow
(\mathop{\rm pr}\nolimits_{{\Bbb A}_{S}^{1}})_{\ast}
\mathop{\rm pr}\nolimits_{X}^{\ast}\rho^{\ast}K\leqno{(\ast )}
$$
puisque $\mathop{\rm pr}\nolimits_{{\Bbb A}_{S}^{1}}\circ\tau
=\mathop{\rm pr}\nolimits_{{\Bbb A}_{S}^{1}}$.

D'apr\`{e}s la formule de K\"{u}nneth, le complexe $(\mathop{\rm pr}
\nolimits_{{\Bbb A}_{S}^{1}})_{\ast} \mathop{\rm pr}
\nolimits_{X}^{\ast}\rho^{\ast}K$ est l'image r\'{e}ciproque du
complexe $f_{\ast}\rho^{\ast}K$ par la projection canonique ${\Bbb
A}_{S}^{1}\rightarrow S$ et la formation de la fl\`{e}che $(\ast )$
commute \`{a} la restriction \`{a} toute section $S\rightarrow {\Bbb
A}_{S}^{1}$ de cette m\^{e}me projection canonique. Il ne reste plus
qu'\`{a} remarquer que la restriction \`{a} la section $1:S\rightarrow
{\Bbb A}_{S}^{1}$ de $(\ast )$ est l'identit\'{e} alors que celle
\`{a} la section $0:S\rightarrow
{\Bbb A}_{S}^{1}$ est nulle pour conclure que
$f_{\ast}\rho^{\ast}K=(0)$.
\hfill\hfill$\square$
\vskip 3mm

\rem D\'{e}monstration de la proposition $5.4$
\endrem
Notons encore $j$ l'inclusion de $V^{\circ}$ dans $V$ et
$\rho^{\circ}:V^{\circ}\rightarrow {\Bbb P}(V)$ le morphisme
quotient. La fl\`{e}che d'oubli des supports de l'\'{e}nonc\'{e} est
un isomorphisme si et seulement si la fl\`{e}che d'oubli des supports
$$
j_{!}\rho^{\ast}K\rightarrow j_{\ast}\rho^{\ast}K
$$
est un isomorphisme dans $D_{{\rm c}}^{{\rm b}}(V,\overline{{\Bbb
Q}}_{\ell})$.  Le lemme d'homotopie pr\'{e}c\'{e}dent appliqu\'{e}
\`{a} $j_{\ast}K\in {\rm ob}D_{{\rm c}}^{{\rm b}}({\cal V},
\overline{{\Bbb Q}}_{\ell})$ assure que cette derni\`{e}re fl\`{e}che
d'oubli des supports est un isomorphisme si et seulement si
$\pi_{\ast}^{\circ}\rho^{\ast}K=\pi_{\ast}j_{\ast}\rho^{\ast}K=(0)$
o\`{u} $\pi^{\circ}=\pi\circ j:V^{\circ}\rightarrow S$.

On a
$$
\pi_{\ast}^{\circ}\rho^{\ast}K=(\overline{\pi}^{\,\circ})_{\ast}(K
\otimes\rho_{\ast}\overline{{\Bbb Q}}_{\ell})
$$
d'apr\`{e}s la formule de K\"{u}nneth, puisque $\rho$ est un ${\Bbb
G}_{{\rm m},S}$-torseur, d'o\`{u} un triangle distingu\'{e}
$$
(\overline{\pi}^{\,\circ})_{\ast}K\rightarrow
\pi_{\ast}^{\circ}\rho^{\ast}K\rightarrow
(\overline{\pi}^{\,\circ})_{\ast}K[-1](-1)\rightarrow .
$$
Par suite, $\pi_{\ast}^{\circ}\rho^{\ast}K=(0)$ si et seulement si la
fl\`{e}che de co-bord
$$
(\overline{\pi}^{\,\circ})_{\ast}K[-1](-1)\rightarrow
(\overline{\pi}^{\,\circ})_{\ast}K[1]
$$
est un isomorphisme, c'est-\`{a}-dire si et seulement si
$(\overline{\pi}^{\,\circ})_{\ast}K=(0)$ puisque ce dernier complexe
est born\'{e}, d'o\`{u} l'assertion (i) de la proposition.

Passons \`{a} la d\'{e}monstration de l'assertion (ii). On a
$$\eqalign{
(\overline{\pi}^{\vee\circ})_{\ast}\mathop{\rm Rad}\nolimits_{{\Bbb
P}(V)/S}(K)&=(\overline{\pi}^{\vee\circ})_{\ast}(q^{\vee})_{\ast}
q^{\ast}K[r-2]\cr
&=(\overline{\pi}^{\circ})_{\ast}q_{\ast}q^{\ast}K[r-2]\cr
&=(\overline{\pi}^{\circ})_{\ast}(K\otimes q_{\ast}\overline{{\Bbb
Q}}_{\ell}[r-2])\cr
&=\bigoplus_{i=0}^{r-2}(\overline{\pi}^{\circ})_{\ast}K[r-2-2i](-i),\cr}
$$
et donc $(\overline{\pi}^{\vee\circ})_{\ast}\mathop{\rm
Rad}\nolimits_{{\Bbb P}(V)/S}(K)=(0)$ si et seulement si
$(\overline{\pi}^{\circ})_{\ast}K=(0)$.

Si $(\overline{\pi}^{\circ})_{\ast}K=(0)$, on a
$$
(j^{\vee})^{\ast}\mathop{\rm Four}\nolimits_{{\cal V}/S}(j_{!}K)
\cong \mathop{\rm Rad}\nolimits_{{\Bbb
P}(V)/S}(K)(-1)
$$
d'apr\`{e}s la proposition 1.6, et donc
$$
j^{\ast}\mathop{\rm Four}\nolimits_{{\cal V}^{\vee}/S}((j^{\vee})_{!}
(j^{\vee})^{\ast}\mathop{\rm Four}\nolimits_{{\cal V}/S}(j_{!}K))\cong
\mathop{\rm Rad}\nolimits_{{\Bbb P}(V^{\vee})/S}(\mathop{\rm
Rad}\nolimits_{{\Bbb P}(V)/S}(K))(-2)
$$
d'apr\`{e}s ce qui pr\'{e}c\`{e}de et encore la proposition 1.6.  Or
on a le triangle distingu\'{e}
$$\displaylines{
\qquad j^{\ast}\mathop{\rm Four}\nolimits_{{\cal V}^{\vee}/S}
((j^{\vee})_{!} (j^{\vee})^{\ast} \mathop{\rm Four}\nolimits_{{\cal
V}/S}(j_{!}K))\rightarrow j^{\ast}\mathop{\rm Four}\nolimits_{{\cal
V}^{\vee}/S}(\mathop{\rm Four}\nolimits_{{\cal
V}/S}(j_{!}K))
\hfill\cr\hfill
\rightarrow j^{\ast}\mathop{\rm Four}\nolimits_{{\cal V}^{\vee}/S}
((i^{\vee})_{\ast} (i^{\vee})^{\ast}\mathop{\rm Four}\nolimits_{{\cal
V}/S}(j_{!}K))\rightarrow ,\qquad}
$$
et on a
$$
j^{\ast}\mathop{\rm Four}\nolimits_{{\cal
V}^{\vee}/S}(\mathop{\rm Four}\nolimits_{{\cal
V}/S}(j_{!}K))\cong K(-r)
$$
d'apr\`{e}s le th\'{e}or\`{e}me 3.1, et
$$
j^{\ast}\mathop{\rm Four}\nolimits_{{\cal V}^{\vee}/S}
((i^{\vee})_{\ast} (i^{\vee})^{\ast}\mathop{\rm Four}\nolimits_{{\cal
V}/S}(j_{!}K))=j^{\ast}\pi^{\ast}(i^{\vee})^{\ast}\mathop{\rm
Four}\nolimits_{{\cal
V}/S}(j_{!}K)[r]=j^{\ast}\pi^{\ast}\pi_{!}j_{!}K[2r]
$$
d'apr\`{e}s des cas particuliers du lemme 1.7 et la proposition 1.8.
Comme le morphisme d\'{e}duit de $\pi\circ j:{\Bbb P}(V)\rightarrow
B({\Bbb G}_{{\rm m},S})$ par le changement de base $S\rightarrow
B({\Bbb G}_{{\rm m},S})$ est le morphisme structural
$\pi^{\circ}:V^{\circ}\rightarrow S$ qui se factorise en le morphisme
quotient $\rho^{\circ} :V^{\circ}\rightarrow {\Bbb P}(V)$ et
$\overline{\pi}^{\circ}:{\Bbb P}(V)\rightarrow S$, il r\'{e}sulte de
l'hypoth\`{e}se $(\overline{\pi}^{\circ})_{\ast}K=(0)$ que
$\pi_{!}j_{!}K=(0)$, et donc que
$$
\mathop{\rm Rad}\nolimits_{{\Bbb P}(V^{\vee})/S}(\mathop{\rm
Rad}\nolimits_{{\Bbb P}(V)/S}(K))\cong K(2-r),
$$
et l'assertion (ii) est d\'{e}montr\'{e}e.
\hfill\hfill$\square$
\vskip 3mm

\thm LEMME 5.6
\enonce
{\rm (i)} La sous cat\'{e}gorie strictement pleine ${\cal S}$ de
$\mathop{\rm Perv}({\Bbb P}(V),\overline{{\Bbb Q}}_{\ell})$
d\'{e}finie comme l'image essentielle du foncteur exact
$$
(\overline{\pi}^{\circ})^{\ast}(\cdot )[r-1]:\mathop{\rm Perv}
(S,\overline{{\Bbb Q}}_{\ell})\rightarrow \mathop{\rm Perv}({\Bbb
P}(V),\overline{{\Bbb Q}}_{\ell}),
$$
est une sous-cat\'{e}gorie de Serre, c'est-\`{a}-dire est stable par
sous-quotients et par extensions.

En particulier, tout $\overline{{\Bbb Q}}_{\ell}$-faisceau pervers $K$
sur ${\Bbb P}(V)$ admet un plus grand sous-faisceau pervers {\rm
(}resp.  faisceau pervers quotient{\rm )} qui provient de $S$ et qui
n'est autre que
$$
(\overline{\pi}^{\circ})^{\ast}({}^{{\rm p}}{\cal
H}^{1-r} (\overline{\pi}^{\circ})_{\ast}K)[r-1]\hookrightarrow K
$$
$$
K\twoheadrightarrow (\overline{\pi}^{\circ})^{!}({}^{{\rm p}}{\cal H}^{r-1}
(\overline{\pi}^{\circ})_{\ast}K)[1-r]=
(\overline{\pi}^{\circ})^{\ast}({}^{{\rm p}}{\cal H}^{r-1}
(\overline{\pi}^{\circ})_{\ast}K)[r-1](r-1)\,)\leqno{{\rm (}resp. }
$$
o\`{u} l'injection {\rm (}resp.  la surjection{\rm )} est induite par
la fl\`{e}che d'adjonction $(\overline{\pi}^{\circ})^{\ast}
(\overline{\pi}^{\circ})_{\ast}\rightarrow \mathop{\rm id}$ {\rm
(}resp.  $\mathop{\rm id}\rightarrow (\overline{\pi}^{\circ})^{!}
(\overline{\pi}^{\circ})_{\ast}${\rm )}.

\decale{\rm (ii)} Notons
$$
{\cal P}({\Bbb P}(V),\overline{{\Bbb Q}}_{\ell})=\mathop{\rm Perv}
({\Bbb P}(V),\overline{{\Bbb Q}}_{\ell})/{\cal S}
$$
la cat\'{e}gorie ab\'{e}lienne quotient au sens de Serre {\rm (}cf.
{\rm [10]} Chapitre {\rm III)}, obtenue en inversant dans la
cat\'{e}gorie ab\'{e}lienne $\mathop{\rm Perv}({\Bbb
P}(V),\overline{{\Bbb Q}}_{\ell})$ les morphismes dont les noyau et
conoyau sont dans ${\cal S}$.  Alors les objets simples de ${\cal P}({\Bbb
P}(V), \overline{{\Bbb Q}}_{\ell})$ sont exactement les objets simples
de $\mathop{\rm Perv}({\Bbb P}(V),\overline{{\Bbb Q}}_{\ell})$ qui ne
sont pas dans ${\cal S}$.

\decale{\rm (iii)} Un $\overline{{\Bbb Q}}_{\ell}$-faisceau pervers
$K$ sur ${\Bbb P}(V)$ est dans ${\cal S}$ si et seulement si
$(\rho^{\circ})^{\ast}K[1]$ est isomorphe dans $\mathop{\rm Perv}
(V^{\circ},\overline{{\Bbb Q}}_{\ell})$ \`{a}
$(\rho^{\circ})^{\ast}K'[1]$ pour un objet $K'$ de ${\cal S}$,
c'est-\`{a}-dire si et seulement si $(\rho^{\circ})^{\ast}K[1]$ est
isomorphe \`{a} $(\pi^{\circ})^{\ast}M'[r]$ pour un objet $M'$ de
$\mathop{\rm Perv}(S,\overline{{\Bbb Q}}_{\ell})$, et on a l'inclusion
$$
(\pi^{\circ})^{\ast}(\mathop{\rm ob}\mathop{\rm Perv}(B({\Bbb G}_{{\rm
m},S}),\overline{{\Bbb Q}}_{\ell}))[r]\subset \mathop{\rm ob}{\cal S}
$$
o\`{u} on a not\'{e} encore par $\pi^{\circ}:{\Bbb P}(V)\rightarrow
B({\Bbb G}_{{\rm m},S})$ le morphisme induit par le morphisme
structural $\pi^{\circ}:V^{\circ}\rightarrow S$.
\endthm

\rem D\'{e}monstration
\endrem
Mis \`{a} part l'assertion de stabilit\'{e} par extensions, la partie
(i) du lemme est d\'{e}montr\'{e}e dans la section 4.2.6 de [6].

Pour tout couple $(M_{1},M_{2})$ de $\overline{{\Bbb
Q}}_{\ell}$-faisceaux pervers sur $S$, le foncteur
$(\overline{\pi}^{\circ})^{\ast}(\cdot )[r-1]$ exact induit un
isomorphisme
$$
\mathop{\rm Ext}\nolimits^{1}(M_{1},M_{2})\buildrel\sim\over
\longrightarrow\mathop{\rm Ext}\nolimits^{1}
((\overline{\pi}^{\circ})^{\ast}M_{1}[r-1],
(\overline{\pi}^{\circ})^{\ast}M_{2}[r-1])
$$
puisque l'on a
$$\eqalign{
\mathop{\rm Ext}\nolimits^{1}((\overline{\pi}^{\circ})^{\ast}M_{1}[r-1],
(\overline{\pi}^{\circ})^{\ast}M_{2}[r-1])&=
\mathop{\rm Ext}\nolimits^{1}((\overline{\pi}^{\circ})^{\ast}M_{1},
(\overline{\pi}^{\circ})^{\ast}M_{2})\cr
&=\mathop{\rm Ext}\nolimits^{1}(M_{1},(\overline{\pi}^{\circ})_{\ast}
(\overline{\pi}^{\circ})^{\ast}M_{2})\cr
&=\mathop{\rm Ext}\nolimits^{1}(M_{1}M_{2}\otimes
(\overline{\pi}^{\circ})_{\ast}\overline{{\Bbb Q}}_{\ell})\cr
&=\mathop{\rm Ext}\nolimits^{1}(M_{1},\bigoplus_{i=0}^{r-1}M_{2}[-2i](-i))\cr
&=\bigoplus_{i=0}^{r-1}\mathop{\rm Ext}\nolimits^{1-2i}(M_{1},
M_{2}(-i))\cr
&=\mathop{\rm Ext}\nolimits^{1}(M_{1},M_{2}).\cr}
$$
Par suite, ${\cal S}$ est stable par extensions , ce qui ach\`{e}ve la
preuve de la partie (i) du lemme.

La partie (ii) est formelle: elle vaut pour toute cat\'{e}gorie
quotient au sens de Serre d'une cat\'{e}gorie ab\'{e}lienne
noeth\'{e}rienne et artinienne.

Enfin, la premi\`{e}re assertion de la partie (iii) r\'{e}sulte de
l'isomorphisme \'{e}vident
$$
K\cong {}^{{\rm p}}{\cal H}^{-1}(\rho_{\ast}^{\circ}
(\rho^{\circ})^{\ast}K[1]),
$$
et la seconde assertion r\'{e}sulte de la premi\`{e}re puisque le
morphisme d\'{e}duit par le changement de base $S \rightarrow B({\Bbb
G}_{{\rm m},S})$ du morphisme $\pi^{\circ}:{\Bbb P}(V)\rightarrow
B({\Bbb G}_{{\rm m},S})$ est le morphisme structural
$\pi^{\circ}:V^{\circ}\rightarrow S$ et se factorise donc en
$\rho^{\circ}:V^{\circ}\rightarrow {\Bbb P}(V)$ et
$\overline{\pi}^{\circ}:{\Bbb P}(V)\rightarrow S$.
\hfill\hfill$\square$
\vskip 3mm

On dira que les objets de ${\cal S}$ sont les $\overline{{\Bbb
Q}}_{\ell}$-faisceaux pervers sur ${\Bbb P}(V)$ {\it qui proviennent
de} $S$.

Bien entendu, on peut remplacer $V$ par $V^{\vee}$ dans ce qui
pr\'{e}c\`{e}de et on obtient la sous-cat\'{e}gorie \'{e}paisse et
stable par extensions ${\cal S}^{\vee}$ de $\mathop{\rm Perv}({\Bbb
P}(V^{\vee}),\overline{{\Bbb Q}}_{\ell})$ dont les objets sont ceux
qui proviennent de $S$, et la cat\'{e}gorie ab\'{e}lienne quotient
au sens de Serre
$$
{\cal P}({\Bbb P}(V^{\vee}),\overline{{\Bbb Q}}_{\ell})=\mathop{\rm
Perv}({\Bbb P}(V^{\vee}),\overline{{\Bbb Q}}_{\ell})/{\cal S}^{\vee}.
$$

\thm PROPOSITION 5.7 (Brylinski [2])
\enonce
{\rm (i)} Si $K$ est un $\overline{{\Bbb Q}}_{\ell}$-faisceau pervers
sur ${\Bbb P}(V)$, alors, pour chaque entier $i\not=0$, le $i$-\`{e}me
faisceau de cohomologie perverse de $\mathop{\rm Rad} \nolimits_{{\Bbb
P}(V)/S}(K)$ provient de $S$, plus pr\'{e}cis\'{e}ment est
canoniquement isomorphe \`{a}
$$
{}^{{\rm p}}{\cal H}^{i}\mathop{\rm Rad}\nolimits_{{\Bbb P}(V)/S}(K)
\cong\openup\jot\cases{(\overline{\pi}^{\vee\circ})^{\ast}({}^{{\rm p}}
{\cal H}^{i-1}(\overline{\pi}^{\circ})_{\ast}K)[r-1]& si $i<0$,\cr
(\overline{\pi}^{\vee\circ})^{\ast}({}^{{\rm p}}{\cal H}^{i+1}
(\overline{\pi}^{\circ})_{\ast}K(1))[r-1]& si $i>0$.\cr}
$$
En particulier, $\mathop{\rm Rad}\nolimits_{{\Bbb P}(V)/S}(K)$ est
pervers si $(\overline{\pi}^{\circ})_{\ast} K\in {}^{{\rm p}}D_{{\rm
c}}^{[-1,1]}(S,\overline{{\Bbb Q}}_{\ell})$.

\decale{\rm (ii)} Le foncteur ${}^{{\rm p}}{\cal H}^{0} \mathop{\rm
Rad}\nolimits_{{\Bbb P}(V)/S}:\mathop{\rm Perv} ({\Bbb P}(V),
\overline{{\Bbb Q}}_{\ell})\rightarrow \mathop{\rm Perv} ({\Bbb
P}(V^{\vee}),\overline{{\Bbb Q}}_{\ell})$ envoie la sous-cat\'{e}gorie
${\cal S}$ dans la sous-cat\'{e}gorie ${\cal S}^{\vee}$ et induit donc
un foncteur not\'{e} ${\cal R}_{{\Bbb P}(V)/S}$ de la cat\'{e}gorie
ab\'{e}lienne ${\cal P}({\Bbb P}(V), \overline{{\Bbb Q}}_{\ell})$ dans
la cat\'{e}gorie ab\'{e}lienne ${\cal P}({\Bbb P}(V^{\vee}),
\overline{{\Bbb Q}}_{\ell})$.  Ce foncteur ${\cal R}_{{\Bbb P}(V)/S}$
est exact et est une \'{e}quivalence de cat\'{e}gories ab\'{e}liennes
de quasi-inverse ${\cal R}_{{\Bbb P}(V^{\vee})/S}(r-2)$.
\endthm

\rem D\'{e}monstration
\endrem
Si $K\in\mathop{\rm ob}\mathop{\rm Perv}({\Bbb P}(V),\overline{{\Bbb
Q}}_{\ell})$, le complexe $(j^{\vee})^{\ast}\mathop{\rm Four}
\nolimits_{{\cal V}/S}(j_{!}K)$ est dans ${}^{{\rm p}}D_{{\rm
c}}^{\leq 0}({\Bbb P}(V^{\circ}), \overline{{\Bbb Q}}_{\ell})$ puisque
$j_{!}$ est $t$-exact \`{a} droite et que $\mathop{\rm Four}
\nolimits_{{\cal V}/S}$ est $t$-exact d'apr\`{e}s le th\'{e}or\`{e}me
4.2. Le triangle distingu\'{e} de la proposition 1.6 donne donc les
isomorphismes cherch\'{e}s pour $i>0$.

Dualement, on a le triangle distingu\'{e}
$$
\mathop{\rm Rad}\nolimits_{{\Bbb P}(V)/S}(K)\rightarrow (j^{\vee})^{\ast}
\mathop{\rm Four}\nolimits_{{\cal V}/S}(j_{\ast}K)(1)\rightarrow
(\overline{\pi}^{\vee\circ})^{\ast}(\overline{\pi}^{\circ})_{\ast}K[r-1]
\rightarrow
$$
et le complexe $(j^{\vee})^{\ast}\mathop{\rm Four}\nolimits_{{\cal
V}/S}(j_{\ast}K)$ est dans ${}^{{\rm p}}D_{{\rm c}}^{\geq 0}({\Bbb
P}(V^{\circ}), \overline{{\Bbb Q}}_{\ell})$, d'o\`{u} les
isomorphismes cherch\'{e}s pour $i<0$ et l'assertion (i) est
d\'{e}montr\'{e}e.

Le calcul
$$\eqalign{
\mathop{\rm Rad}\nolimits_{{\Bbb P}(V)/S}
((\overline{\pi}^{\circ})^{\ast}M[r-1])&=
(q^{\vee})_{\ast}q^{\ast}(\overline{\pi}^{\circ})^{\ast}M[2r-3]\cr
&=(q^{\vee})_{\ast}(q^{\vee})^{\ast}(\overline{\pi}^{\vee\circ})^{\ast 
}M[2r-3]\cr
&=(\overline{\pi}^{\vee\circ})^{\ast}M\otimes
(q^{\vee})_{\ast}\overline{{\Bbb Q}}_{\ell}[2r-3]\cr
&=\bigoplus_{i=0}^{r-2}(\overline{\pi}^{\vee\circ})^{\ast}M[2r-3-2i](-i)}
$$
implique la formule
$$
{}^{{\rm p}}{\cal H}^{0}\mathop{\rm Rad}\nolimits_{{\Bbb P}(V)/S}
((\overline{\pi}^{\circ})^{\ast}M[r-1])=
\openup\jot\cases{(\overline{\pi}^{\vee\circ})^{\ast}M[r-1](-r')
& si $r-2=2r'$ est pair,\cr
(0) & si $r$ est impair,\cr}
$$
pour tout $\overline{{\Bbb Q}}_{\ell}$-faisceau pervers $M$ sur $S$,
et donc l'inclusion
$$
{}^{{\rm p}}{\cal H}^{0}\mathop{\rm Rad}\nolimits_{{\Bbb P}(V)/S}(
\mathop{\rm ob}{\cal S})\subset \mathop{\rm ob}{\cal S}^{\vee}.
$$

Compte tenu du triangle distingu\'{e} de la proposition 1.6 et de la
partie (iii) du lemme 5.7, on a aussi l'inclusion
$$
{}^{{\rm p}}{\cal H}^{0}(j^{\vee})^{\ast}\mathop{\rm Four}
\nolimits_{V/S}(j_{!}(\mathop{\rm ob}{\cal S}))\subset \mathop{\rm ob}
{\cal S}^{\vee}
$$
et le foncteur induit par ${}^{{\rm p}}{\cal H}^{0}\circ
j^{\vee})^{\ast}\circ\mathop{\rm Four} \nolimits_{V/S}\circ j_{!}$ de
${\cal P}({\Bbb P}(V),\overline{{\Bbb Q}}_{\ell})$ dans ${\cal
P}({\Bbb P}(V^{\vee}), \overline{{\Bbb Q}}_{\ell})$ n'est autre que
${\cal R}_{{\Bbb P}(V)/S}$. On conclut en raisonnant comme dans la
d\'{e}monstration de la partie (ii) de la proposition 5.4.
\hfill\hfill$\square$
\vskip 3mm

\thm COROLLAIRE 5.8 (Brylinski [2])
\enonce
{\rm (i)} Si $K$ est un $\overline{{\Bbb Q}}_{\ell}$-faisceau
pervers {\rm irr\'{e}ductible} sur ${\Bbb P}(V)$ qui ne provient pas
de $S$, le faisceau pervers $L={}^{{\rm p}}{\cal H}^{0}\mathop{\rm
Rad}\nolimits_{{\Bbb P}(V)/S}(K)$ sur ${\Bbb P}(V)$ admet une
filtration \`{a} trois crans
$$
(0)=L_{0}\subset L_{1}\subset L_{2}\subset L_{3}=L
$$
telle que $L_{1}$ et $L_{3}/L_{2}$ proviennent de $S$ et que
$L_{2}/L_{1}$ soit un $\overline{{\Bbb Q}}_{\ell}$-faisceau pervers
irr\'{e}ductible sur ${\Bbb P}(V^{\vee})$.

\decale{\rm (ii)} Si $K$ est un $\overline{{\Bbb Q}}_{\ell}$-faisceau
pervers {\rm irr\'{e}ductible} sur ${\Bbb P}(V)$ tel que
$(\overline{\pi}^{\circ})_{\ast}K=(0)$, alors $\mathop{\rm Rad}
\nolimits_{{\Bbb P}(V)/S}(K)$ est un $\overline{{\Bbb
Q}}_{\ell}$-faisceau pervers irr\'{e}ductible sur ${\Bbb P}
(V^{\vee})$.
\endthm

\rem D\'{e}monstration
\endrem
L'assertion (i) r\'{e}sulte de l'assertion (ii) de la proposition 5.7
puisque celle-ci assure que toute suite de Jordan-H\"{o}lder de $L$
admet un et un seul sous-quotient qui ne provient pas de $S$.

On sait d\'{e}j\`{a} que $\mathop{\rm Rad}\nolimits_{{\Bbb P}(V)/S}
(K)$ est pervers d'apr\`{e}s l'assertion (i) de la proposition 5.7.
Consid\'{e}rons une filtration
$$
(0)=L_{0}\subset L_{1}\subset L_{2}\subset
L_{3}=L=\mathop{\rm Rad}\nolimits_{{\Bbb P}(V)/S}(K)
$$
comme dans l'assertion (i) d\'{e}j\`{a} d\'{e}montr\'{e}e.  Si $L_{1}=
(\overline{\pi}^{\vee\circ})^{\ast}M_{1}[r-1]$ pour un
$\overline{{\Bbb Q}}_{\ell}$-faisceau pervers $M_{1}$ sur $S$, on a
$$
\mathop{\rm Hom}(L_{1},L)=\mathop{\rm Hom}
(M_{1}[r-1],(\overline{\pi}^{\vee\circ})_{\ast}L)
$$
dans $D_{{\rm c}}^{{\rm b}}(S,\overline{{\Bbb Q}}_{\ell})$.  Par
suite, on a $\mathop{\rm Hom}(L_{1},L)=(0)$, et donc n\'{e}cessairement
$L_{1}=(0)$, puisque $(\overline{\pi}^{\vee\circ})_{\ast}L=(0)$ (cf.
l'assertion (ii) de la proposition 5.4).  Comme les hypoth\`{e}ses sur
$K$ sont auto-duales, on voit de m\^{e}me que $L/L_{2}=(0)$, d'o\`{u}
le corollaire.
\hfill\hfill$\square$
\vskip 5mm

{\bf  Bibliographie}
\vskip 5mm

\newtoks\ref	\newtoks\auteur
\newtoks\titre	\newtoks\editeur
\newtoks\annee	\newtoks\revue
\newtoks\tome	\newtoks\pages
\newtoks\reste	\newtoks\autre

\def\bibitem#1{\parindent=8pt\itemitem{#1}\parindent=12pt}

\def\livre{\bibitem{[\the\ref]}%
\the\auteur ~-- {\sl\the\titre}, \the\editeur, ({\the\annee}).
\smallskip\smallskip\filbreak}

\def\article{\bibitem{[\the\ref]}%
\the\auteur ~-- \the\titre, {\sl\the\revue} {\the\tome},
({\the\annee}), \the\pages.\smallskip\filbreak}

\def\autre{\bibitem{[\the\ref]}%
\the\auteur ~-- \the\reste.\smallskip\filbreak}

\ref={1}
\auteur={E. {\pc FRENKEL}, D. {\pc GAITSGORY}, K. {\pc VILONEN}}
\titre={On the geometric Langlands conjecture}
\revue={J. Amer. Math. Soc.}
\tome={15}
\annee={2002}
\pages={367-417}
\article

\ref={2}
\auteur={J.-L. {\pc BRYLINSKI}}
\reste={Transformations canoniques, Dualit\'{e} projective,
Th\'{e}orie de Lefschetz, Transformations de Fourier et Sommes
trigonom\'{e}triques, dans {\it G\'{e}om\'{e}trie et Analyse
Microlocales, Ast\'{e}risque} 140-141, (1986), 3-134}
\autre

\ref={3}
\auteur={G. {\pc LAUMON}}
\reste={Travaux de Frenkel, Gaitsgory et Vilonen sur la correspondance
de Drinfeld-Langlands, S\'{e}minaire Bourbaki, juin 2002, 54\`{e}me
ann\'{e}e, 2001-2002, n\raise4pt\hbox{o} 906}
\autre

\ref={4}
\auteur={T. {\pc EKEDAHL}}
\reste={On the Adic Formalism, dans {\it The Grothendieck Festschrift,
Volume II, Progess in Mathematics} 87, {\it Birkh\"{a}user},
(1990), 197-218}
\autre

\ref={5}
\auteur={G. {\pc LAUMON}, L. {\pc MORET}-{\pc BAILLY}}
\titre={Champs alg\'{e}briques}
\editeur={Springer-Verlag}
\annee={1999}
\livre

\ref={6}
\auteur={A. {\pc BEILINSON}, J. {\pc BERNSTEIN}, P. {\pc DELIGNE}}
\reste={Faisceaux pervers, {\it Ast\'{e}risque} 100, (1982)}
\autre

\ref={7}
\auteur={O. {\pc GABBER}}
\reste={Notes on some $t$-structures, manuscrit, d\'{e}cembre
2000}
\autre

\ref={8}
\auteur={G. {\pc LAUMON}}
\titre={Transformation de Fourier, constantes d'\'{e}quations
fonctionnelles et conjecture de Weil}
\revue={Publ. Math. de l'I.H.\'{E}.S.}
\tome={65}
\annee={1987}
\pages={131-210}
\article

\ref={9}
\auteur={N.M. {\pc KATZ}, G. {\pc LAUMON}}
\titre={Transformation de Fourier et majoration de sommes
exponentielles}
\revue={Publ. Math. de l'I.H.\'{E}.S.}
\tome={62}
\annee={1986}
\pages={145-202}
\article

\ref={10}
\auteur={P. {\pc GABRIEL}}
\titre={Des cat\'{e}gories ab\'{e}liennes}
\revue={Bull. Soc. math. France}
\tome={90}
\annee={1962}
\pages={323-448}
\article
\bye